\documentclass[12pt]{article}
\def\llex {\leq_{_{ d-LEX}}}
\def\sllex{<_{_{ d-LEX}}}
\def\ring{k}
\def\contr{\neg}
\def\transp{^T}
\def\ka{K\"{a}hler~}
\def\contr{\rightharpoonup}
\def\st{\mathbf{^sT}}
\def\hka{Hyper\"{a}hler}

\def\rnum{\mathbf{R}}
\def\cnum{\mathbf{C}}
\def\znum{\mathbf{Z}}
\def\sl{\mathbf{sl}(s+1,\rnum)}
\def\slc{\mathbf{sl}(s+1,\cnum)}

\def\gs{\mathbf{sl}(s+1,\mathbf{C})}

\def\ls{\mathbf{l}_{s}}

\def\tpm1{\bigotimes^{p-1}\,\! _{\ring}}

\def\wpm1{\bigwedge^{p-1}\,\!\!_{\ring}}

\def\primq{\mathbf{Prim}^{q}}

\def\prim{\mathbf{Prim}}
\def\ei{{E^{i}}}
\def\eik{{E_{k}^{i}}}
\def\ii{{I^{i}}}
\def\iik{{I_{k}^{i}}}
\def\ks{{\mathbf{k}_{s}}}

\def\me{\mathbf{g}}
\def\eme{\bar{\mathbf{g}}}
\def\qd  {\begin{picture}(8,8)
  
 \put(0,0){\line (0,1){7}}
 \put(0,0){\line (1,0){7}}
 \put(7,7){\line (0,-1){7}}
 \put(7,7){\line (-1,0){7}}
  
 \end{picture}}
 \def \ep {\hfill \qd}
\newtheorem{teo}{Theorem}[section]

\newtheorem{cor}[teo]{Corollary}
\newtheorem{exe}[teo]{Example}
\newtheorem{lem}[teo]{Lemma}
\newtheorem{pro}[teo]{Proposition}
\newtheorem{dfn}[teo]{Definition}

\newtheorem{rmk}[teo]{Remark}
\begin{document}
\title{Polysymplectic spaces, s-\ka manifolds and Lagrangian Fibrations}
\author{Michele Grassi}
\date{June 21, 2000}
\maketitle

\begin{abstract}
In this paper we begin the study of polysymplectic manifolds, and of their relationship with PDE's. This notion provides a generalization of symplectic manifolds which is very well suited for the geometric study of PDE's with values in a smooth manifold. Some of the standard tools of analytical mechanics, such as the Legendre transformation and Hamilton's equations, are shown to generalize to this new setting. There is a strong link with lagrangian fibrations, which can be used to build Polysymplectic manifolds.\\
We then provide the definition and some basic properties of s-K\"{a}hler and almost s-K\"{a}hler 
manifolds. These are a generalization of the usual notion of K\"{a}hler and almost K\"{a}hler
manifold, and they reduce to them for $s = 1$. The basic 
properties of K\"{a}hler manifolds, and their Hodge theory, can be generalized to s-K\"{a}hler manifolds, 
with some modifications.  The most interesting examples come from semi-flat special lagrangian fibrations of Calabi-Yau manifolds.
\end{abstract}
\tableofcontents

\textbf{Introduction}\\
In this paper we begin the study of polysymplectic manifolds, and of their relationship with PDE's. A (non degenerate) polysymplectic manifold  of rank $n$ is a smooth manifold, together with $s$ closed differential $2$-forms on it such that, if we indicate the data with $(\mathcal{M},\omega_1,...,\omega_s)$, for any point $p\in \mathcal{M}$ there exist near $p$ coordinates $x_1,..,x_n,y^1_1,...,y^s_n$ providing a description $\omega_j = \sum_i dx_1\wedge dy^j_i$. From the definition, it is clear that for $s = 1$ a polysymplectic manifold is just a symplectic manifold. This definition is not the one which we will give in the paper, but it has the advantage of being very quick and explicit.  The notion just introduced  provides a generalization of the notion of symplectic manifold which is very well suited, as we will show in this paper, to the study of PDE's with values in a smooth manifold.
While some results from the theory of symplectic manifolds generalize to the polysymplectic setting, and others do not, it is almost always the case that the proofs in the $s > 1$ case are different in spirit from the ones when $s = 1$. This is true for example for Theorem ~\ref{teo-multisympnorm}, which generalizes Darboux's Theorem, where one has to face same questions of integrability which were not present in the symplectic setting.\\
When dealing with PDE's with values in a smooth manifold $M$ and $s$ independent variables, the natural ambient space is $TM\times_M\cdots\times_MTM$ ($s$ times), which we will for simplicity indicate with $\st(M)$. We will also use the notation $\st^*(M)$ for the analogous construction obtained starting from $T^*(M)$. It turns out that for PDE's which come from a Lagrangian (i.e. a smooth function $\mathcal{L}$ on $\st(M)$), the usual machinery of analytical mechanics can be generalized rather effortlessly. In particular, from any Lagrangian we obtain a Legendre transformation, which maps $\st(M)$ to $\st^*(M)$.
A Lagrangian is non degenerate if and only if  a suitable polysymplectic structure on $\st(M)$ associated to $\mathcal{L}$ is non degenerate, and in this case the Legendre transformation associated to it is a local diffeomorphism. The picture becomes even more compelling after one notices that on $\st^*(M)$ there is a canonical non degenerate polysymplectic structure, and that the Legendre transformation maps the polysymplectic structure on $\st(M)$ associated canonically to any non degenerate lagrangian to the canonical polysymplectic structure on $\st^*(M)$. The Euler-Lagrange equations for a Lagrangian $L$ are then mapped  to a set of equations resembling in a striking way the Hamilton equations. Of course, a polysymplectic manifold needs not have a non degenerate Poisson  structure (it could be odd dimensional, for example!); however, there is a way to put in a canonical way on any polysymplectic manifold a generalization of a non degenerate Poisson structure, involving not one but $s$ "brackets". The generalization of the Hamilton equations mentioned earlier can be expressed in terms of the canonical $s$-Poisson structure on $\st^*(M)$, in a way that resembles the classical one (see Theorem ~\ref{gener-ham}).\\
There is an apparently different approach to the geometry of PDE's. For a guide to that approach, see for example ~\cite{G} or consult the bibliography of ~\cite{G2} for further references. It would be interesting to investigate the possible links between these two viewpoints.

In the second part of the paper we introduce the notions of {\em almost s-\ka}  and of {\em s-\ka} manifold. s-\ka manifolds are a generalization of \ka manifolds, to which they reduce when $s=1$. A smooth manifold $M$ of dimension $n(s+1)$ together with a 
Riemannian metric $\mathbf{g}$ and $2$-forms 
$\omega_{1},...,\omega_{s}$ is {\em s-\ka} if the data satisfy the 
following property:
for each point  $p\in M$ there exist an open neighborhood $\mathcal{U}$ 
of $p$ and a system of coordinates 
$x_{i},y^{j}_{i}$,$i~=~1,...,n$, $j~=~1,...,s$
on $\mathcal{U}$ such that with these coordinates:\\
1)~$\omega_{j}~=~\sum_{i=1}^ndx_{i}\wedge dy^{j}_{i} ~~ ~\forall j \in\{1,...,s\}$,\\
2)~ $\mathbf{g}_{\alpha,\beta}~=~\delta^\alpha_\beta~+~\mathbf{O}(2)$\\
They are "rigid" objects, with an extremely rich set of properties, which generalize those of \ka manifolds. In the last part of the paper we explore some of them.
A part of the structure of a s-\ka manifold, namely the forms $\omega_1,...,\omega_s$ described above, determines a polysymplectic structure. Polysymplectic manifolds are much "softer" objects, with no local moduli, but with a rich global geometry. In between these two notions, there is the notion of an {\em almost s-\ka manifold}, which seems to share some of the structure of s-\ka manifolds, without being too rigid. Any  polysymplectic manifold admits a Riemannian metric for which it is almost s-\ka, in exactly the same way as any symplectic manifold admits a metric compatible with the symplectic structure, which makes it almost complex.
In this paper we show that one can put an almost 2-\ka structure on interesting compact 3-manifolds. At this stage of the theory, we don't see any general obstruction on putting a polysymplectic structure on any 3-manifold. Polysymplectic threefolds need not be orientable either. Because a polysymplectic structure has non trivial global invariants, this might be of some help in the classification of 3-manifolds.
The most interesting examples by far of almost s-\ka manifolds that we build in this paper are those that come from special lagrangian fibrations of Calabi-Yau manifolds. The constructions of the current paper work in the so-called {\em semi-flat} case described in ~\cite{SYZ}, but in a forthcoming paper we will show that one can deal also with the general situation. We refrained form including a treatment of this topic, and of the relationship of all this with mirror symmetry (in the spirit of ~\cite{SYZ}), due to the already excessive size of the present paper.\\
Going back to s-\ka manifolds, we prove that they enjoy many of the properties of \ka manifolds, and in particular we start the development of their Hodge theory. We prove that there is an analogue to the Hodge identities, and to the Lefschetz decomposition, namely a representation of the Lie algebra $\mathbf{sl}(s+1,\rnum)$ on the cohomology of any compact oriented s-\ka manifold. We also show that there is an analogue to the hard Lefschetz theorem. The representation mentioned above is induced by one on differential forms, and reduces to the standard one when we are in the \ka case, namely for $s=1$. Although many facts concerning (almost) s-\ka and polysymplectic manifolds resemble similar properties of (almost) \ka and symplectic manifolds, the proofs in the non classical situations are almost always very different from the classical ones, and exploit new properties. 
This paper is unfortunately very rich of explicit computations, because we had to establish directly many basic properties of the objects we introduced, for lack of a reference. We hope that this work will allow us to write more readable papers on this subject in the future, as we will be able to refer to this computations without having to reproduce them. It should be noted that this paper  has been written over a period of more than three years, with many interruptions; this implies that there will be many more misprints that we would like. We hope that the reader will be willing to let us know of any misprint/mistake that he or she might find.\\
We now come to a more detailed description of the contents of the paper:\\
In the first section we introduce polysymplectic vector spaces, and we prove a normal form theorem for them. We then prove a Lemma which describes the qualitative structure of any change between sets of coordinates in which the polysymplectic structure is in normal form (standard coordinates). In the same Lemma we introduce some elements of the second wedge power of the space, which will be used later to define s-Poisson manifolds.\\
In the second section we first introduce polysymplectic manifolds, and after  giving the central example $\st^*(M)$, we prove the local normal form Theorem for polysymplectic manifolds, which generalizes Darboux's Theorem. In the last part of the section we introduce the canonical s-Poisson structure on $\st^*(M)$ compatible with a given (non degenerate) polysymplectic structure, and prove its existence.\\
In the third section we provide some examples of polysymplectic manifolds, the most notable being the one that starts from lagrangian fibrations.\\
In the fourth section we introduce the notion of compatibility between a Riemannian metric and a polysymplectic structure. We then proceed to the proof of the theorem which shows that for any given polysymplectic structure the space of metrics compatible with it is non empty and contractible.\\
In the fifth section we introduce first order Lagrangians (also called $s$-Lagrangians) on $\st(M)$, and define the canonical polysymplectic structure on $\st(M)$ associated to a non degenerate Lagrangian. We then define the "energy" $\mathcal{H}_\mathcal{L}$ associated to a strongly non degenerate Lagrangian $\mathcal{L}$, and the Legendre transform associated to any Lagrangian. We finally prove an equivalent formulation of the Euler-Lagrange equations associated to a Lagrangian, and use it (in the non degenerate case) to translate them via the Legendre transform in a generalization of the Hamilton equations on $\st^*(M)$.\\
In the sixth section we introduce almost s-\ka manifolds, and analyze their relationship with polysymplectic ones, via the notion of "compatibility" of a metric with a polysymplectic structure.
We give some examples of the above notions, for manifolds with assigned topology (but not compact), and for compact three-manifolds. The main theorem of this section is the last one, describing the strong relationship between special lagrangian fibrations of Calabi-Yau manifolds and almost $s$-\ka manifolds.\\
In section ~\ref{sec:ska} we introduce s-\ka manifolds, and we show that they can be seen as almost s-\ka manifolds with an integrability condition added.\\
In section ~\ref{sec:exeska} we give some examples of s-\ka manifolds. Apart from \ka ones, the simplest nontrivial examples that we give are manifolds which are diffeomorphic to (real) tori. Although this establishes that there are compact s-\ka manifolds for all s, it is not satisfactory, and we expect to be able to provide more interesting examples in the future. Using special lagrangian fibrations of Calabi-Yau manifolds which enjoy some (rather strong) flatness conditions, we are able to build other examples of $s$-\ka manifolds. This construction will be analyzed in much more detail in a forthcoming paper, in which we will investigate the connection of this circle of ideas (and the ones described in the following  sections) with the approach of ~\cite{SYZ} to mirror symmetry.\\ 
In section ~\ref{sec:lefop} we introduce a generalization of the Lefschetz operators on s-\ka manifolds, and analyze their commutation relations with their adjoints with respect to the metric. The explicit computation of the commutation relations allows us to show that these operators, together with their adjoints, generate a Lie algebra isomorphic (canonically) to $\mathbf{sl}(s+1,\rnum)$.\\
In section ~\ref{sec:hodgeska} we prove a generalization of the Hodge identities, which allow us to show that the action of $\mathbf{sl}(s+1,\rnum)$ introduced in the previous section induces an action on the cohomology of any compact oriented s-\ka manifold. We introduce primitive forms, and show that there is an analogue of the Lefschetz decomposition.\\
In the eleventh section we prove that not only there is a Lefschetz decomposition, but one can also prove a generalization of the hard Lefschetz theorem.\\
In the last section we draw some conclusions, concerning both polysymplectic manifolds and s-\ka manifolds.\\
\textbf{Notations}\\
For a vector space $V$ and an element $\alpha\in \bigwedge^*V^*$, we will indicate with $C(\alpha)$ the smallest subspace of $V^*$ such that $\alpha \in \bigwedge^*(C(\alpha))$. For an $\alpha\in \bigwedge^*V^*$, we indicate with $A(\alpha)\subset V$ the subspace orthogonal to $\alpha$ with respect to contraction.
The remaining notations will be either standard or explicitely introduced in the paper. A note on terminology: we found out that the name "s-symplectic" is widely used to indicate an object which does not have anything to do with what we introduce here. This was the reason for our use of the (less appealing) terminology "polysymplectic". We have recently found out that also this name has been used, although much less diffusely. We hope that this overlap will not cause any problems. 
\section{Polysymplectic vector spaces}
\label{sec:polyvectors}
In this section we define polysymplectic vector spaces, and study some of their properties. The most important fact is Theorem ~\ref{multsympvect-nf}. Lemma ~\ref{injectivityssymp} will be used in the next section to define the generalization of Poisson manifolds. In this section the base field $\ring$ will be assumed to be the Real numbers, but many results would continue to hold on any field of characteristic different from $2$.
\begin{dfn}
Let $V$ be a vector space over $\ring$, and let $\omega_{1},...,\omega_{s}$ be forms in $\bigwedge^{2}(V^{*})$. We say that the forms induce a {\em polysymplectic} structure on $V$ of rank $n$ if the following three conditions hold:\\
$1)$ The forms
\[\omega_{1}^{\wedge i_{1}}\wedge\cdots\wedge\omega_{s}^{\wedge i_{s}}\]
where we vary the s-tuple of non-negative integers 
$( i_{1},..., i_{s})$ subject to the condition $ i_{1}+\cdots+ i_{s}~=~n$ are all independent.\\
$2)$ The forms
\[\omega_{1}^{\wedge i_{1}}\wedge\cdots\wedge\omega_{s}^{\wedge i_{s}}\]
where we vary the s-tuple of non-negative integers 
$( i_{1},..., i_{s})$ subject to the condition $ i_{1}+\cdots+ i_{s}~=~n+1$ are all $0$. \\
$3)$ (If $s>1$) For all $j\in\{1,...,s\}$, $dim\left(C(\omega_j)\bigcap \sum_{k\not= j}C(\omega_k)\right) \leq n$.\\
The polysymplectic structure is said {\em non-degenerate} if $dim(V) = s(n+1)$.
\end{dfn}
\begin{exe}
A vector space with a nondegenerate antisymmetric bilinear form (i.e. a symplectic vector space) is polysymplectic
\end{exe}
\begin{rmk}
Let $(V,\omega_{1},...,\omega_{s})$ be a polysymplectic vector space of rank $n$. Then for any number $t$ such that $1\leq t < s$, $(V,\omega_{1},...,\omega_{t})$ is a polysymplectic vector space of the same rank.
\end{rmk}
We now show that a polysymplectic vector space can be put in a normal form, in the same way as a symplectic vector space has a basis in which the non-degenerate two-form has a canonical expression.
\begin{dfn}
The forms
$\omega_{1},\ldots,\omega_{s}~\in~\bigwedge^{2}V^{*}$
are said to be in {\em polysymplectic normal form} with respect to a basis of $V$ if the basis is of the form $e_{1},\ldots,e_{n},e^{1}_{1},\ldots,e^{s}_{n},z_1,...,z_{dim(V)-(s+1)n}$
and the forms can be expressed in terms of this basis as
\[\omega_{i}~=~\sum_{j=1}^n e_{j}^{*}\wedge {e^{i}_{j}}^{*}\]
We call such a basis {\em polysymplectic} or {\em standard} for  
$\left(V, \omega_{1},\ldots,\omega_{s}\right)$
\end{dfn}
We will use the following standard fact from symplectic linear algebra:
\begin{pro} Let $V$ be a vector space (over a field of characteristic different from $2$), $\omega\in\bigwedge^2V^*$, $\omega_n\not=0$, $\omega_{n+1}=0$, and let $e_1,...,e_t,f_1,..,f_r~\in V$ be vectors such that $<e_1,...,e_t,f_1,...,f_r>\bigcap \omega^\perp = (0)$, and for the bilinear form $B_\omega$ associated to $\omega$ we have for all $i_1,i_2<t,j_1,j_2<r$ that $B_\omega(e_{i_1},e_{i_2}) = 0$, $B_\omega(f_{j_1},f_{j_2}) = 0$, $B_\omega(e_{i_1},f_{j_1}) = \delta_{i_1 j_1}$, then there are vectors $e_{t+1},...,e_n$, $f_{r+1},...,f_n$, $z_1,..z_{dim(V)-2n} \in V$ such that for all $i_1,i_2\leq n$, ~$j_1,j_2\leq n$ and $k_1,k_2\leq  dim(V)-2n$, we have $B_\omega(e_{i_1},e_{i_2}) = 0$, $B_\omega(f_{j_1},f_{j_2}) = 0$, $B_\omega(e_{i_1},f_{j_1}) = \delta_{i_1 j_1}$, $B_\omega(e_{i_1},z_{k_1}) = B_\omega(f_{j_1},z_{k_1}) = B_\omega(z_{k_1},z_{k_2}) = 0$.
\end{pro}
\begin{teo}
\label{multsympvect-nf}
Let $V$ be a vector space, and
$\omega_{1},\ldots,\omega_{s}~\in~\bigwedge^{2}V^{*}$.
Then the following are equivalent:\\
1) The forms $\omega_{1},\ldots,\omega_{s}$ determine a polysymplectic structure on $V$.\\
2) $V$ has a polysymplectic basis with respect to the forms $\omega_1,...,\omega_s$.
\end{teo}
{\it Proof of $2)\Rightarrow 1)$}\\
Suppose that the forms 
$\omega_{1},\ldots,\omega_{s}~\in~\bigwedge^{2}V^{*}$
can be expressed, in terms of a basis
$e_{1},\ldots,e_{n},e^{1}_{1},\ldots,e^{s}_{n}~\in~V$
as
$\omega_{i}~=~\sum_{j}e_{j}^{*}\wedge {e^{i}_{j}}^{*}$ (i.e. the basis is polysymplectic with respect to the forms $\omega_1,...,\omega_s$).
Then we have that:\\
$1)$ For any $(h_{1},\ldots,h_{s})$ such that $h_{i}~\geq~0$ and
$\sum_{i}h_{i}~=~n$ 
the forms
\[\left\{w_{1}^{\wedge h_{1}}\wedge\cdots\wedge w_{s}^{\wedge 
h_{s}}~|~h_{i}~\geq~0~ and
~\sum_{i}h_{i}~=~n\right\}\]
are independent (and in particular all different from zero).\\
Indeed, suppose $\sum_{|H|=n}\alpha_H\omega^H~=~0$, where we have used the multi index notation $\omega^H = w_{1}^{\wedge h_{1}}\wedge\cdots\wedge w_{s}^{\wedge 
h_{s}}$ and $|(h_1,..,h_s)| = h_1+\cdots+h_s$. Then, if $H_0$ is the lowest (lexicographically) multi index such that $\alpha_H\not= 0$, the form $\omega^{H_0}$ contains the "monomial" 
\[(e_1^*\wedge\cdots\wedge e_n^*)\wedge ({e^{1*}_1}\wedge\cdots \wedge {e^{1*}_{h_1}})\wedge\cdots\wedge ({e^{s*}_1}\wedge\cdots\wedge {e^{s*}_{h_s}})\]
with coefficient $\pm 1$. It is immediate to check that this monomial appears only once in the expression for $\omega^{H_0}$, and never appears in any $\omega^{H}$ for $H > H_0$. Therefore, it must be $\alpha_{H_0} = 0$, contradiction.\\
$2)$
$w_{1}^{\wedge h_{1}}\wedge\cdots\wedge w_{s}^{\wedge h_{s}}$
whenever $\sum_{i}h_{i}~>~n$\\
Indeed, any "monomial" in $\omega^H$ for $|H|>n$ must contain at least $|H|$ forms from the set $\{e_1^*,...,e_n^*\}$. Any such wedge product must therefore be zero.\\

\textit{Proof of $1) \Rightarrow 2)$}
\begin{lem}
Let $(V,\omega_1,\omega_2)$ be a polysymplectic vector space of rank $n$. Then $dim_\ring\left(C(\omega_1)\bigcap C(\omega_2)\right) = n$
\end{lem}
\textit{Proof}
Given the forms $\omega_{1}$ and $\omega_{2}$, take any decomposition
\[V~=~V_{1}\oplus V_{2}\oplus V_{3}\oplus V_{4}\]
where
\[V_{1}~=~A(\omega_{1})\cap A(\omega_{2}),~ V_{1}\oplus V_{2}~=~A(\omega_{1}),~
V_{1}\oplus V_{3}~=~A(\omega_{2})\]
We then have that 
\[V^{*}~=~ V_{1}^{*}\oplus V_{2}^{*}\oplus V_{3}^{*}\oplus V_{4}^{*}\]
and 
\[C(\omega_{1})~=~ V_{3}^{*}\oplus V_{4}^{*},~ ~
C(\omega_{2})~=~ V_{2}^{*}\oplus V_{4}^{*}\]
It follows that
\[ C(\omega_{1})\cap C(\omega_{2})~=~ V_{4}^{*}\]
W now prove two facts concerning the above decomposition:\\
1)
$\omega_{1}|_{V_{4}}~=~0,\omega_{2}|_{V_{4}}~=~0$\\
Consider the first statement (the proof of the second one is the same).
If the restriction of $\omega_{1}$ to $V_{4}$ is not $0$, it follows that we can build a standard basis $e_1,...,e_n,f_1,..,f_n,z_1,...z_{n(s-1)}$ for it which has $e_1,f_1\in V_4$, and where $V_1\oplus V_2$ is the span of $z_1,...z_{n(s-1)}$. It is easy to check that 
this contradicts the dimension  condition on the mixed $n^{th}$ wedge products of $\omega_1$ with $\omega_2$.\\
2) $dim(V_{4})~\geq~n$\\
We already know that $\omega_{1}|_{V_{4}}~=~0,\omega_{2}|_{V_{4}}~=~0$. Then, if $dim(V_{4})~<~n$
we can build a standard basis $e_1,...,e_n$, $f_1,..,f_n$, $z_1,...z_{n(s-1)}$ for $\omega_{1}$  where $V_{4}$ is the span of $e_1,...,e_a$, with $a~<~n$, and $V_1\oplus V_2$ is the span of $z_1,...z_{n(s-1)}$. It is easy to check that from this it follows that 
$\omega_{1}^{\wedge n}\wedge\omega_{2}~\not=~0$.\\
The two facts above prove that on one hand
\[dim\left(C(\omega_{1})\cap C(\omega_{2})\right)~\leq~n,\]
while on the other hand 
\[dim\left(C(\omega_{1})\cap C(\omega_{2})\right)~\geq~n\]
and therefore we obtain the thesis.
\ep

We now observe that from the lemma and point $3)$ of the definition of a polysymplectic structure, it follows that \[dim\left(C(\omega_{r})\cap\sum_{i\not=r} C(\omega_{i})\right)~=~n~ for ~all ~r\in\{1,...,s\}\]
This, together with the fact that $dim\left(C(\omega_r)\cap C(\omega_j)\right) = n$ for $j\not= r$, implies that $dim\bigcap_j C(\omega_j) = n$.
We may therefore take a basis $\phi_1,...,\phi_n$ of $\bigcap_jC(\omega_j)$. You then have that this basis can be completed with $\psi_i^j$ 
($i=1,...,n,j=1,..,s$) such that $\omega_j~=~\sum_i \phi_i\wedge\psi^j_i$,
 and from the dimension assumption it follows that the $\phi_i,\psi^j_i$ are independent. Complete this basis with $\zeta_1,...,\zeta_{dim(V)-n(s+1)}$, and take now $e_i,e^j_i,z_k$ to be the dual basis to $\phi_i,\psi^j_i,\zeta_k$. It is clear by construction that $e_i,e^j_i$ is a standard basis for $\omega_1,...,\omega_s$.
\ep

\begin{cor}
There exists a representation of the permutation group over $s$ elements on $V$, which induces the permutation of the forms $\omega_{1},...,\omega_{s}$ on $\bigwedge^{2}V$.
There exists a representation of the permutation group over $n$ elements on $V$, which leaves the $s$ forms fixed (and is not trivial). 
The two above representations commute.
\end{cor}
\begin{lem}
\label{injectivityssymp}
Let $(V,\omega_{1},...,\omega_{s})$ be a  vector space with a non-degenerate polysymplectic structure. Let $\sigma_j~:~V\to V^*$ be the operator which contracts a vector with the two-form $\omega_j$. We then have that:\\
1) $\bigoplus_j\sigma_j:~V \to \bigoplus_jV^*$ is injective.\\
2)  If $s>1$ and $e_{1},\ldots,e_{n}$, $e^{1}_{1},\ldots,e^{s}_{n}$ 
and $f_{1},\ldots,f_{n}$, $f^{1}_{1},\ldots,f^{s}_{n}$ are two 
polysymplectic bases, there exist an invertible $n\times n$ matrix 
$\mathbf{\theta} = (\theta^m_i)$ and a tensor $\eta = \eta^{m j}_{i}$ such that 
\[f_i = \sum_m\theta^m_ie_m~+~\sum_{m,j}\eta^{m j}_{i}e^{j}_{m},~ ~f^j_i = \sum_m 
\left(\mathbf{\theta}^{-1}\right)^i_m e^j_m~ ~for~all~i,j\]
and with $\sum_m\theta_i^m\eta^{mj}_k =  \sum_m\theta_k^m\eta^{mj}_i$.
If moreover the two bases are orthogonal with respect to some 
(positive definite non degenerate) metric, the tensor $\eta$ must vanish identically.
\end{lem}
\textit{Proof}
Pick a standard polysymplectic basis $e_{1},\ldots,e_{n}$, $e^{1}_{1},\ldots,e^{s}_{n}$, so that $\omega_j = \sum_{i=1}^n e_i^*\wedge e^{j*}_i$. \\
1) Assume that $\sigma_j(v) = 0$ for all $j$. If we express $v$ using  the standard basis, we see that we can write $v = v_j + w_j$, with $v_j$ involving only $e_{1},\ldots,e_{s},e^{j}_{1},\ldots,e^{j}_{n}$, and $w_j$ involving only the other basis vectors. It is clear that $\sigma_j(w_j) = 0$, and therefore from $\sigma_j(v) = 0$ we get $\sigma_j(v_j) = 0$. At this point we conclude that $v_j = 0$, from the standard fact that a symplectic form is non degenerate. Repeating this argument for all $j\in\{1,..,s\}$, we see that $\forall j ~v_j = 0$, and hence $v = 0$.\\
2) Assume that $f_{1},\ldots,f_{n}$, $f^{1}_{1},\ldots,f^{s}_{n}$ is another polysymplectic basis for $\omega_1,...,\omega_s$. Because the span of $e_1^j,...,e^j_n$ is $\bigcap_{k\not= j}\omega_k^\perp$, we have that \[<e_1^j,...,e^j_n> = <f_1^j,...,f^j_n>~ ~ ~for ~all ~j\in\{1,...,s\}\]
There are therefore matrices $\mathbf{\theta} = (\theta_i^l)$ $\mathbf{\theta}^j = (\theta_i^{j,m})$, and a tensor $\eta = (\eta^{m j}_{i})$  such that
$f_i = \sum_l\theta_i^le_l + \sum_{m,j}\eta^{m j}_{i}e^k_m$ and $f_i^j = \sum_m\theta_i^{j,m}e^j_m$. Because $B_{\omega_j}(f_i,f_m^j) = \delta_{im}$, we see that $\sum_l\theta^l_i\theta^{j,l}_m = \delta_{im}$, or in other words for all $j$ we have that  $\mathbf{\theta}^j = (\mathbf{\theta}^{-1})^{\transp}$. \\
We know that $B_{\omega_j}(f_i,f_k) = 0$ for all $i,k,j$. Writing up what this means in terms of $\theta$ and $\eta$, wee see that 
\[0 = B_{\omega_j}(f_i,f_k) = \sum_m\left(\theta_i^m\eta^{mj}_k - \theta_k^m\eta^{mj}_i\right)\]
which proves the statement. 
The proof of the last part of point $2)$ is immediate.\\
\ep
\begin{cor}
Let $(V,\omega_{1},...,\omega_{s})$ be a  vector space with a non-degenerate polysymplectic structure, and $s > 1$. Then there is a canonically determined subspace $C^V$ of $V^*$ of dimension $n$, characterized by the property of being the span of $e_{1}^*,\ldots,e_{n}^*$, for some (and therefore any) polysymplectic basis $e_{1},\ldots,e_{n}$, $e^{1}_{1},\ldots,e^{s}_{n}$. Dualizing, there is a canonically determined subspace $\left(C^V\right)^\perp$ of $V$ of dimension $ns$, characterized by the property of being the span of $e^{1}_{1},\ldots,e^{s}_{n}$ for some (and therefore any) polysymplectic basis $e_{1},\ldots,e_{n}$, $e^{1}_{1},\ldots,e^{s}_{n}$.
\end{cor}

\section{Polysymplectic and s-Poisson manifolds}

\begin{dfn}
\label{dfn-polysymplectic}
A manifold $M$ of dimension $n(s+1)$ together with $s$ closed degree $2$ differential forms $\omega_{1},...,\omega_{s}$ is said to be polysymplectic if \\
1) For each point $p$ in $M$ 
\[(T_{p}M,{(\omega_{1}})_{p},...,({\omega_{s}})_{p})\]
is a polysymplectic vector space;\\
2) The distribution of subspaces $C^M\subset T^*M$ defined for $p\in M$ as
\[C^M_p = \left\{\alpha\in T^*_pM~:~\alpha\wedge(\omega_1^n)_p = \cdots = \alpha\wedge(\omega_s^n)_p = 0 \right\}\]
is generated (locally) by closed smooth $1$-forms.\\
We say that the polysymplectic manifold is {\em non-degenerate} if the polysymplectic vector spaces $(T_{p}M,{\omega_{1}}_{p},...,{\omega_{s}}_{p})$ are all non-degenerate.
\end{dfn}
In the following, unless otherwise stated, all the polysymplectic manifolds will be assumed to be non-degenerate.
\begin{exe}
A symplectic manifold is polysymplectic (with $s = 1$)
\end{exe}
\begin{exe}
\label{execotangent}
Let $M$ be a smooth manifold, and let $\mathbf{T}^{*}M$ indicate the cotangent bundle of $M$. If 
\[\st^*(M)~:=~\mathbf{T}^{*}M~ ~\times_{_{M}} ~~\cdots~ ~\times_{_{M}} ~~\mathbf{T}^{*}M~(s~times),\]
\[\st(M)~:=~\mathbf{T}M~ ~\times_{_{M}} ~~\cdots~ ~\times_{_{M}} ~~\mathbf{T}M~(s~times),\]
$\pi_{i}~:~\st^*(M)\to \mathbf{T}^{*}M$
is the projection on the $i^{th}$ factor, and $\omega$ the canonical symplectic form on 
$\mathbf{T}^{*}M$, let
$\omega_{i}~:=~ \pi_i^{*}\omega$.
We have then that $(\st^*(M),\omega_{1},...,\omega_{s})$ is polysymplectic.
\end{exe}
{\it Proof}\\
The forms $\omega_i$ are closed, because they are pull-back of closed forms. Moreover, point by point they induce a polysymplectic structure on $T\left(\st^*(M)\right)$, because if we choose coordinates $x_1,...,x_n$ on $\mathcal{U}\subset M$ around $p$, we can pick coordinates $x_1,...,x_n,y^1_1,...,y^s_n$ on $\st^*(\mathcal{U})\subset\st^*(M)$, such that in these coordinates \[\omega_j~=~\sum_{i=1}^n dx_i\wedge dy^j_i\] 
It is then clear that at every point the frame dual to the coframe $dx_1,...,dx_n$, $dy^1_1,...,dy^s_n$ is a standard polysymplectic basis for $\omega_1,...,\omega_s$. To conclude, it is enough to observe that in the open set where the coordinates are defined 
\[C^{\st^*(M)} ~= ~< dx_1,...,dx_n >\]
and is therefore locally generated by closed forms, as required by the definition.
\ep
\begin{dfn}
Let $M$ be a smooth manifold. If  $x_1,...,x_n$ is a system of coordinates on $\mathcal{U}\subset M$,  the coordinates $x_1,...,x_n,y^{1}_{1},...,y^{s}_{n}$ on $\st^*(\mathcal{U})\subset \st^*(M)$, correspond to the frame $\left(\frac{\partial }{\partial x_{1}}\right)_1,...,\left(\frac{\partial }{\partial x_{n}}\right)_1,...,\left(\frac{\partial }{\partial x_{1}}\right)_s,...,\left(\frac{\partial }{\partial x_{n}}\right)_s$ for the vector bundle $\st(\mathcal{U})$ over $\mathcal{U}$.
\end{dfn}

\begin{teo}[Polysymplectic normal form]
\label{teo-multisympnorm}
Let $(X,\omega_1,...,\omega_s)$ be a smooth polysymplectic manifold and $p\in X$. Assume given elements $\phi_1,...,\phi_n,\psi^1_1,...,\psi^n_s$ of $T^*_pX$ such that for all $j=1,..,s$ one has $(\omega_j)_p=\sum_i\phi_i\wedge\psi^i_j$.Then there are a neighborhood $\mathcal{U}\subset X$ of $p\in X$, a neighborhood $\mathcal{V}\subset \mathbf{R}^{dim(X)}$ of $0\in \mathbf{R}^{dim(X)}$ and an isomorphism of polysymplectic manifolds
\[\phi~:~(\mathcal{U},\omega_1,...,\omega_s)~\rightarrow~
\left(\mathcal{V},\sum_i dx_i\wedge dy^1_i,...,\sum_i dx_i\wedge dy^s_i\right)\]
where we indicated the coordinates on $\mathbf{R}^{dim(X)}$ with $x_1,...,x_n,y^1_1,...,y^s_n$. With this notation, one can also assume $(dx_i)_p=\phi$, $(dy^i_j)_p=\psi^i_j$.
\end{teo} 
{\it Proof}\\
If $V$ is a vector space, given an element of $\alpha\in\bigwedge^*(V)$ we indicate with $C(\alpha)$ the smallest subspace $W\subset V$ such that $\alpha\in \bigwedge^*W$. Similarly, for a differential form $\alpha$ we define $C(\alpha)$ to be the smallest distribution of subspaces $D\subset\Omega^1$ such that $\alpha\in \bigwedge^*D$. A priori, the $C(\omega_j)$ are only "generalized Pfaffian systems", as defined for example in ~\cite[Page 382]{LM}. From Darboux's Reduction Theorem, in the form stated for example in ~\cite[Bryant, Page 103]{FU}, we see that $C(\omega_j)$ is a vector bundle (of rank $2n$) for any $j=1,...,s$, with local coframes given by closed $1$-forms. We clearly have that $C^X=\bigcap_j C(\omega_j)$. Then $C^X$ is a constant rank distribution of subspaces of $T^*X$, which by the definition of a polysymplectic structure is locally generated by closed forms. From the constant rank property, we may assume that there are (locally) $n$ functions $x_1,...,x_n$ such that $dx_1,...,dx_n$ are independent, and for all $q$ in the open set considered
\[< (dx_1)_p,\dots,(dx_1)_p > = C^X\]
By acting if necessary with a constant transformation matrix we can assume that $\forall i~ (dx_i)_p=\phi_i$.\\
Fix now an index $j\in\{1,...,s\}$. From Darboux's reduction theorem, we can find coordinates $z_1,...,z_d$ such that $\omega_j$ is expressed only in terms of $z_{d-2n+1},...,z_d$, and such that $\frac{\partial}{\partial z_k}$ is in $C(\omega_j)^\perp$ for $k=1,..,d-2n$ (and therefore one has also $<dz_{d-2n+1},...,dz_d>=C(\omega_j)$). From their definition, it follows that $\frac{\partial x_i}{\partial z_k}=0$ for all $i$, and for $k=1,..,d-2n$. Therefore, we can apply the theorem of Carath\'{e}odory-Jacobi-Lie (see ~\cite[Page 136]{LM}) to conclude that there are functions $y_i^j$ (depending only on the $z_{d-2n+1},...,z_d$)  such that $dy^j_i\in C(\omega_j)$ and $\omega_j~=~
\sum_i dx_i\wedge dy^j_i$. Because $<dx_1,...,dx_n,dy^1_j,...,dy^n_j>=C(\omega_j)$, by an invertible linear transformation inside $C(\omega_j)$ (with constant coefficients) leaving all the $dx_i$ fixed  we can also assume that $dy^i_j=\psi^i_j$. After repeating the procedure for all $j$, we end up with functions $x_1,..,x_n,y^1_1,...,y^s_n$ near $p\in X$.
The $x_1,...,x_n,y^1_1,...,y^s_n$ form a system of coordinates because the $dx_1,...,dx_n,dy^1_1,...,dy^s_n$ are independent forms.
\ep
\begin{cor}
Let $M$ be a smooth manifold, and $\omega_1,...,\omega_2$ be smooth $2$-forms on it. The following are then equivalent:\\
1) $(M,\omega_1,...,\omega_s)$ is a polysymplectic manifold.\\
2) For all $p\in M$ there are coordinates $x_1,..,x_n,y^1_1,...,y^s_n$ near $p$ such that 
\[\forall j~ ~\omega_j = \sum_i dx_i\wedge dy^j_i\]
\end{cor}
Note that the second condition of the previous corollary was given as a definition of a polysymplectic manifold in the introduction.
\begin{cor}
If $s+1$ is even, any polysymplectic manifold $(\mathcal{M},\omega_1,...,\omega_2)$ is orientable
\end{cor}
\textit{Proof} Recall that we assume that the polysymplectic structure is non degenerate. The statement is then a direct consequence of the previous theorem and of Lemma ~\ref{injectivityssymp}, part $2)$.
\ep
\begin{teo}
\label{simpsdiff2}
Let $M$ be a smooth manifold, and let $\omega_1,...,\omega_s$ be 
smooth $2$-forms on it, with $s\not= 2$. Let $n\geq 1$ be an integer. 
The following are equivalent:\\
1) The forms $\omega_1,...,\omega_s$ determine a non degenerate 
polysymplectic structure on $M$\\
2) The forms $\omega_1,...,\omega_s$ are closed, and induce a non 
degenerate polysymplectic structure on $T_pM$ for all $p\in M$
\end{teo}
As is shown in Example ~\ref{contr-s2}, the theorem fails for $s=2$.\\
\textit{Proof}\\
We clearly have only to prove that $2)\Rightarrow 1)$, and this only 
for $s \geq 3$, because for $s=1$ it is immediate. By looking at the 
definition of polysymplectic manifold, we see that we have to prove 
that 
the distribution of subspaces $C^M\subset T^*M$ defined for $p\in M$ 
as
\[C^M_p = \left\{\alpha\in T^*_pM~:~\alpha\wedge(\omega_1^n)_p = 
\cdots = \alpha\wedge(\omega_s^n)_p = 0 \right\}\]
is generated (locally) by closed smooth $1$-forms.
From the fact that the forms $\omega_{1},\ldots,\omega_{s}$ induce polysymplectic structures of 
constant rank $n$ on the various $T_pM$, we deduce that $C^M$ is a 
constant rank distribution of subspaces of $T^{*}M$ (of rank $n$).
\begin{lem}
\label{smoothdistr}
In the notation above, $C^M$ is a differentiable distribution, i.e. 
is generated by smooth sections.
\end{lem}
{\it Proof of the lemma} To see that $C^M$ is differentiable, 
trivialize locally $T^*M$, so that it looks like $\mathbf{R}^d\times 
(\mathbf{R}^d)^*$. The generalized distribution $C^M$ is then a 
smooth family of affine subspaces of this vector space. Take a 
subspace $S_v$ of $\mathbf{R}^d\times (\mathbf{R}^d)^*$ of dimension 
$2d-n$, and which intersects $(C^M)_0$ in a given point $v$. Then for 
any $q\in \mathbf{R}^d$ near enough to $0$,
$S_v\cap (C^M)_q$ will be formed by just one point $\sigma_v(q)$. It 
is clear from the construction that $\sigma_v$ is a smooth section of 
$C^M$ near $p$, passing through $v$. If we vary $v$ through a basis 
of $T_pM$, the $\sigma_v$ so obtained generate $C^M$ in a 
neighborhood of $p$. Note that here we used the fact that $C^M$ has 
constant rank. 
\ep

The local sections of $C^M$ generate an (algebraic) ideal inside 
$\Omega^*$, which we indicate with $I^M$. From the theorem of 
Frobenius, it is now enough to show that $dI^M\subset I^M$.
For this, assume that $\alpha\in\Omega^1(M)$ is a section of $C^M$, 
and therefore an element of $I^M$. By definition, 
\[\alpha\wedge(\omega_1^n) = \cdots = \alpha\wedge(\omega_s^n) = 0\]
By differentiating this equations, using the fact that the $\omega_j$ 
are closed, we obtain that
\[(d\alpha)\wedge(\omega_1^n) = \cdots = (d\alpha)\wedge(\omega_s^n) 
= 0\]
Using Theorem ~\ref{multsympvect-nf}, it is easy to prove that any 
$2$-form $\beta$ satisfying $\beta\wedge(\omega_1^n) = \cdots = 
\beta\wedge(\omega_s^n) = 0$ must be a section of $C^M\wedge 
T^*M\subset \bigwedge^2T^*M$. For this we use the fact that $s\geq 3$.
We now put a Riemannian metric on $M$, and apply the Gram-Schmidt process to 
conclude that not only $C^M$ is 
differentiable, but it has also a complementary differentiable 
distribution (its orthogonal distribution with respect to the chosen 
metric). Once that is done, it is clear that if $F$ is such a 
complementary distribution, any smooth section $\beta$ of $C^M\wedge 
T^*M$ can be uniquely written as  a sum 
$\sum_k\alpha_k\wedge\gamma_k$, with the $\alpha_k$ sections of $C^M$ 
and the $\gamma_k$ sections of $T^*M$. \\
Summing up, we have proved that if $\alpha\in I^M\cap \Omega^1$, then 
$d\alpha\in I^M\wedge \Omega^1$. This is clearly enough to show that 
$dI^M\subset I^M$, which proves, via the Frobenius Theorem, that 
$C^M$ is generated locally by closed forms.
\ep

\begin{dfn} 1) Let $M$ be a smooth manifold, and let \[\{~,~\}_j~:~\mathcal{C}^\infty(M)\times \mathcal{C}^\infty(M)\to \mathcal{C}^\infty(M)\] 
be local bilinear skew symmetric maps which are derivations with respect to each factor, for $j\in\{1,...,s\}$. The maps $\{~,~\}_j$ define an {\em s-Poisson structure}  on $M$ if for any $p\in M$ there exist an open neighborhood $\mathcal{U}$ of $p$ in $M$ and local coordinates $x_1,...,x_n,y^1_1,...,y^s_n$ on $\mathcal{U}$ such that 
$\{f,g\}_j = \sum_i \left(\frac{\partial f}{\partial x_i}\frac{\partial g}{\partial y^j_i} - \frac{\partial g}{\partial x_i}\frac{\partial f}{\partial y^j_i}\right)$ on $\mathcal{U}$ for all $j\in\{1,...,s\}$ and all $f,g\in \mathcal{C}^\infty_c(\mathcal{U})$\\
2) Let $(M,\omega_{1},...,\omega_{s})$ be a polysymplectic manifold. A s-Poisson structure 
on $M$ is {\em compatible with the polysymplectic structure} if for any $p\in M$ there exist an open neighborhood $\mathcal{U}$ of $p$ in $M$ and local coordinates $x_1,...,x_n,y^1_1,...,y^s_n$ on $\mathcal{U}$ such that 
\[\{f,g\}_j = \sum_i \left(\frac{\partial f}{\partial x_i}\frac{\partial g}{\partial y^j_i} - \frac{\partial g}{\partial x_i}\frac{\partial f}{\partial y^j_i}\right)~ ~and~ ~\omega_j = \sum_i dx_i\wedge dy^j_i\] 
on $\mathcal{U}$ for all $j\in\{1,...,s\}$ and all $f,g\in \mathcal{C}^\infty_c(\mathcal{U})$   
\end{dfn} 
\begin{teo}
\label{canspoisson}
Let $M$ be a smooth manifold, and let $(\st^*(M),\omega_1,...,\omega_s)$ be the polysymplectic manifold canonically associated to it. There exists then a canonical s-Poisson structure $\{~,~\}_1,....\{~,~\}_s$ on $\st^*(M)$ compatible with the given polysymplectic structure.
\end{teo}
\textit{Proof} Fix any polysymplectic frame $\frac{\partial}{\partial x_1},...,\frac{\partial}{\partial x_n},\frac{\partial}{\partial y_1^1},....,\frac{\partial}{\partial y^s_1} $of $T_p\left(\st^*(M)\right)$, such that $\frac{\partial}{\partial x_1},...,\frac{\partial}{\partial x_n}$ span the tangent space to the zero section of $\st^*(M)$, when seen as a vector bundle over $M$. We can define point by point sections $\alpha_j \in \bigwedge^2T_p\left(\st^*(M)\right)$ as ${\alpha_j}_{|_{T_pM}} = \sum_i \frac{\partial}{\partial x_i}\wedge \frac{\partial}{\partial y_i^j}$. Using Lemma ~\ref{injectivityssymp}, it is easy to see that any other polysymplectic frame $\frac{\partial}{\partial w_1},...,\frac{\partial}{\partial w_n},\frac{\partial}{\partial z_1^1},....,\frac{\partial}{\partial z^s_1} $of $T_p\left(\st^*(M)\right)$, such that $\frac{\partial}{\partial w_1},...,\frac{\partial}{\partial w_n}$ span the tangent space to the zero section of $\st^*(M)$, will be related to the previous one by a linear change of coordinates of the form
\[\frac{\partial}{\partial w_i} = \sum_m\theta^m_i \frac{\partial}{\partial x_m},~ ~ ~\frac{\partial}{\partial z_i^j} = \sum_m(\theta^{-1})^i_m \frac{\partial}{\partial y_m^j}\]
for some invertible $n\times n$ matrix $\theta$. We have therefore that
\[\sum_{i=1}^n \frac{\partial}{\partial w_i} \wedge \frac{\partial}{\partial z_i^j} = 
\sum_{i=1}^n \left(\sum_m\theta^m_i \frac{\partial}{\partial x_m}\right)\wedge \left(\sum_r (\theta^{-1})^i_r \frac{\partial}{\partial y_r^j}\right) =\] \[\sum_{m,r}\left(\sum_i\theta^m_i(\mathbf{\theta}^{-1})^i_r\right)\frac{\partial}{\partial x_m}\wedge \frac{\partial}{\partial y_r^j} =
\sum_{m,r}\delta_{mr} \frac{\partial}{\partial x_m}\wedge \frac{\partial}{\partial y_r^j} = \alpha_j\]
To see that these sections are smooth, using the definition of the canonical polysymplectic structure on $\st^*(M)$, take a system of coordinates $x_1,...,x_n,y_1^1,...,y^s_n$ in a neighborhood $\mathcal{U}$ around $p\in \st^*(M)$, such that for all $j$ we have $\omega_j = \sum_i dx_i\wedge dy^j_i$ on $\mathcal{U}$, and the $x_1,...,x_n$ are coordinates on the zero section of $\st^*(M)$, when seen as a vector bundle over $M$.. It is then clear that ${\alpha_j}_{|_{T_qM}} = \sum_i \frac{\partial}{\partial x_i}\wedge \frac{\partial}{\partial y_i^j}$ for  $q\in\mathcal{U}$, for all $j$. This proves smoothness, and also that if we use the sections $\alpha_j$ to induce bilinear maps on smooth functions via $\{f,g\}_j~:=~<f\otimes g,\alpha_j>$, the $\{~,~\}_j$ satisfy all the properties of an s-Poisson structure compatible with the polysymplectic structure $\omega_{1},...,\omega_{s}$.
\ep
\section{Examples and Lagrangian fibrations}
\label{sec:examplespolysymp}
Before building our first examples, we need the following
\begin{dfn}
A {\em Lagrangian fibration} is a smooth map $f$ from a polysymplectic manifold  $(X,\omega_1,...,\omega_s)$ to a smooth manifold $B$, such that (for $s>1$) for any point $p\in X$ one has $Ker(df)_p = \sum_j(\omega_j^\perp)$. For $s=1$ we simply require that $\omega_{|_{Ker(df)_p}} = 0$.
\end{dfn}
\begin{exe}[Lagrangian Fibrations]
\label{teopolysymp}
Let $(X,\omega^X_1,...,\omega^X_{s_X})$,  $(Y,\omega^Y_1,...,\omega^Y_{s_Y})$ be polysymplectic manifolds of dimensions respectively $(s_X+1)n$ and $(s_Y+1)n$.
Let $f:X\to B$ and $g:Y\to B$ be lagrangian fibrations with smooth fibres, with $dim(B) = n$. Consider the fibred product manifold $M = X\times_B Y \subset X\times Y$, and let $\pi_X$ (resp. $\pi_Y$)  be the projections on the first (resp. second) factor. Let
$\omega_j = \pi_X^*\omega^X_j$ for $j\leq s_X$, and $\omega_j = \pi_Y^*\omega^Y_{j-s_X}$ for $s_X <j\leq s_X + s_Y$. Then $(M,\omega_1,...,\omega_{s_X+s_Y})$ is polysymplectic, and the natural map $F:M\to B$ is a lagrangian fibration with smooth fibres.
\end{exe}
\textit{Proof}\\
We have to verify three things: that the $\omega_j$ are closed, that they induce a polysymplectic structure point by point on $M$, and that the distribution $C^M$ inside $T^*M$ is generated locally by closed $1$-forms.\\
The fact that $M$ is a manifold, and that the $\omega_j$ are closed forms on it, is clear from the definitions.\\
 To verify that the $\omega_j$ induce a polysymplectic structure point by point on $M$, let $\underline{p} = (p,q)$ be a point in $M$. Pick metrics on $T_pX$ and $T_qY$ which are compatible with the polysymplectic structures. If $s_X=1$ (resp. $s_Y = 1$) we also require that the induced almost complex structure in that case swaps $Ker(df)_p$ (resp. $ker(dg)_q$) with its orthogonal complement. 
Pick $n$ independent tangent vectors $u_1,...,u_n$ in $T_{f(p)}B = T_{g(q)}B$, and let $v_1,..., v_n$ (resp. $w_1,...,w_n$) be vectors in $T_pX$ (resp. $T_qX$) such that $df(v_j) = dg(w_j) = u_j$ for all $j\in\{1,...,n\}$. The $v_j,w_j$ are uniquely defined by the further property that we take them to be orthogonal to the fibre. We have therefore, by definition, that the $(v_j,w_j)$ are vectors in $T_{\underline{p}}M$ for $j=1,...,n$.\\
Observe that given a basis  $v_1,...,v_n$ of $Ker(df)_p^\perp$ (resp. $w_1,...,w_n$ of $Ker(dg)_q^\perp$), we can always complete it to a polysymplectic basis $v_1,...,v_n,v^1_1,...,v^{s_X}_n$ for $T_pX$  (resp. $w_1,...,w_n,w^1_1,...,w^{s_Y}_n$ for $T_pY$) using only vectors along the fibre. This is an easy consequence of the structure results on posysymplectic vector spaces. However, it is clear that with the obvious identifications
\[(v_1,w_1),...,(v_n,w_n),(v_1^1,0),...,(v_n^{s_X},0),(w_1^1,0),...,(w_n^{s_Y},0)\]
is a polysymplectic basis for $T_{\underline{p}}M$.
The last thing that we have to verify is that the distribution $C^M$ is locally generated by closed forms. To see this, observe that there is a naturally induced smooth map $F:M\to B$, and that $C^M$ is the distribution of forms orthogonal to $Ker(dF)$; as such it is locally generated by closed forms from Frobenius' theorem (as $Ker(dF)$ is an integrable distribution).
The fact that the map $F$ is a lagrangian fibration with smooth fibres follows from the fact that it is smooth, and from the pointwise description of its differential given above.

\ep
\begin{exe}[Open subsets]
If a manifold $M$ is an open subset of a manifold $(N,\omega_1,...,\omega_s,\me)$ which is polysymplectic (respectively almost s-\ka), then in obtains by restriction a polysymplectic (respectively an almost s-\ka) structure. \\
This applies, for example, when $M$ in a proper open subset of $S^n$, because in that case it can also be considered as an open subset of $D^n$, via stereographic projection.
\end{exe}
\begin{rmk}
Let $M$ be a smooth manifold, and $(T^*(M),\omega)$ its cotangent bundle, with the canonical symplectic form. Then the natural map $T^*(M)\to M$ is a lagrangian fibration, and the multisymplectic structure on $\st^*(M)$ described in Example ~\ref{execotangent} is obtained by iterating ($s$ times) the construction above.
\end{rmk}
Note that by using $s$-cotangent bundles we can construct polysymplectic manifolds with assinged homotopy typle (albeit not compact).

\begin{exe}[Compact polysymplectic manifolds of dimension 3]
\label{2ka-3m} ~\\
Consider the  "solvmanifold" $G_1/\Gamma_1$, defined in ~\cite[Theorem 1.9 page 73 and  Example 2.2 page 77]{TO}. If $\mathbf{g}_1$ is the Lie algebra of $G_1$,
\[\mathbf{g}_1 ~=~<~X,Y,Z;~[X,Y] = kY,~[X,Z] = -kZ,~[Y,Z] = 0~>\]
(with $k\in\rnum$ fixed), pick as  $\omega_1$ and $\omega_2$ the invariant forms induced by the Lie algebra cocycles $X^*\wedge Y^*$ and $X^*\wedge Z^*$ respectively. The forms $\omega_1$ and $\omega_2$ determine a polysymplectic structure on $G_1/\Gamma_1$.
\end{exe}
\textit{Proof}
For completeness, we write a description of $G_1/\Gamma_1$. The reader is advised to go to ~\cite[Theorem 1.9 page 73 and  Example 2.2 page 77]{TO} for details. $G_1 = \rnum \times_\phi \rnum^2$, with 
\[\phi(t) = \left(\begin{array}{cc} e^{kt} & 0 \\
0 & e^{-kt}\end{array}\right)\]
in some basis for $\rnum^2$, with $k\in\rnum$, $e^{kt} + e^{-kt}\not= 2$ an integer, and $\Gamma_1 = \znum \times_\phi \znum^2$.\\
It is immediate to see that the invariant forms $\omega_1$ and $\omega_2$, induced by the Lie algebra cocycles $X^*\wedge Y^*$ and $X^*\wedge Z^*$ respectively, are closed, invariant, and define point by point a polysymplectic structure on $T_x(G_1/\Gamma_1)$ for all $x\in G_1/\Gamma_1$. To see that they induce a polysymplectic structure on $G_1/\Gamma_1$ it is enough to observe, in view of Theorem ~\ref{simpsdiff2}, that $C^{G_1/\Gamma_1}$ is spanned by the invariant closed form induced by the cocycle $X^*$.
\ep
\begin{exe}[A counterexample to Theorem ~\ref{simpsdiff2} for $s=2$]
\label{contr-s2} ~\\
On $SL(2,\rnum)$, take $\omega_1$ and $\omega_2$ to be the invariant forms induced by the cocycles $h^*\wedge e^*$ and $h^*\wedge f^*$ respectively, for some (fixed) Chevalley basis $e,f,g$ of $\mathbf{sl}(2,\rnum)$. Then $\omega_1$ and $\omega_2$ are closed, invariant, and define a polysymplectic structure on $T_x(SL(2,\rnum))$ for all $x\in SL(2,\rnum)$. However, they do not determine a polysymplectic structure on $SL(2,\rnum)$.
\end{exe}
\textit{Proof}\\
It is enough to show that the distribution $C^{SL(2,\rnum)}$ associated to $\omega_1$ and $\omega_2$ 
(see Definition ~\ref{dfn-polysymplectic}) is not generated locally by closed forms. To see this, note first of all that $C^{SL(2,\rnum)}$ is generated by the invariant form associated to the element $h^*\in \mathbf{sl}(2,\rnum)^*$, and that the de Rham differential of this form is the invariant form $\alpha$ associated to the element $e^*\wedge f^*\in \bigwedge^2\mathbf{sl}(2,\rnum)^*$. If $\phi\alpha$ is any other 1-form generating $C^{SL(2,\rnum)}$ on some open set, with $\phi\in C^\infty(SL(2,\rnum))$, we see that $d(\phi\alpha) = d(\phi)\wedge \alpha ~+~ \phi d(\alpha)$. When we evaluate this expression over any point in the open set, we see that to get zero the left hand side and the right hand side must vanish simultaneously,  and therefore in particular $\phi$ must vanish identically. This proves that no form generating (locally) the distribution $C^{SL(2,\rnum)}$ can be closed, and therefore $\omega_1$ and $\omega_2$ do not define a polysymplectic structure.
\ep

\section{Compatible metrics}
\label{sec:compmetric}
In this section we study the relation of compatibility between a metric and a polysymplectic structure. We show that given a polysymplectic structure, there is always a nonempty contractible space of metrics compatible with it.
\begin{dfn}
Let $(M,\omega_{1},...,\omega_{s})$ be a polysymplectic manifold. A 
Riemannian metric $\mathbf{g}$ on $M$ is {\em compatible with the 
polysymplectic structure} if for each point $p$ of $M$ there exists 
an orthonormal basis of $T^*_pM$, 
$dx_{i},dy^{j}_{i}$ (with $i~=~1,...,n,$ $j~=~1,...,s$)
such that 
$\forall j ~(\omega_{j})_{p}~=~\sum_{i}dx_{i}\wedge dy^{j}_{i}$.
\end{dfn}
This section is devoted to the proof of the following Theorem, which is standard for $s =1$, i.e. in the case of symplectic manifolds.
\begin{teo}
Let $(M,\omega_1,...,\omega_s)$ be a (non-degenerate) polysymplectic manifold. The space of Riemannian metrics on $M$ compatible with the polysymplectic structure is non-empty and contractible.
\end{teo}
For the purposes of this proof, we give the following definition. We will use it again when dealing with special lagrangian fibrations.
\begin{dfn}
Let $(V,\omega_1,...,\omega_s)$ be a vector space with a non-degenerate polysymplectic structure, $s > 1$. A Riemannian metric $\mathbf{g}$ on $V$ is {\em block-compatible} with the polysymplectic structure it there exists a polysymplectic basis $e_1,...,e_n,$ $f^1_1,...,f^s_n$ such that for all $i,m,j,k$ (with $j\not= k$)
\[\mathbf{g}(e_i,f^j_m) = \mathbf{g}(f^k_i,f^j_m) = 0\]
\end{dfn}
\begin{lem}
Let $(V,\omega_1,...,\omega_s)$ be a vector space with a non-degenerate polysymplectic structure, $s > 1$, and let $\mathbf{g}_1$ and $\mathbf{g}_2$ be two Riemannian metrics on $V$ block-compatible with the polysymplectic structure, and such that their restrictions to the span of the spaces $\omega_j^\perp$  coincide. If $t\in [0,1]$, then the Riemannian metric $t\mathbf{g}_1 + (1-t)\mathbf{g}_2$ is also block-compatible with the polysymplectic structure.
\end{lem}
\textit{Proof}\\
In view of the block-compatibility of the two metrics with the polysymplectic structure and of Lemma \ref{injectivityssymp}
there are vectors $d_1,...,d_n$, $f_1,...,f_n$, $f^1_1,...,f^s_n$ such that $e_1,...,e_n,f^1_1,...,f^s_n$ and $d_1,...,d_n,f^1_1,...,f^s_n$ are polysymplectic bases, and for $j \not= k$
\[\mathbf{g}_1(e_i,f^j_m) = \mathbf{g}_1(f^k_i,f^j_m) = 0, \mathbf{g}_2(d_i,f^j_m) = \mathbf{g}_2(f^k_i,f^j_m) = 0\]
Moreover,  we can take for all $j$ bases $h^j_1,...,h^j_n$ of the span of $f^j_1,...,f^j_n$, orthonormal with respect to $\mathbf{g}_1$ (and therefore also with respect to $\mathbf{g}_2$). We do not require such bases $h^j_1,...,h^j_n$ to be part of a polysymplectic basis. Such a basis exists because of the hypothesis on the behavior of the two metrics on the span of $f^1_1,...,f^s_n$.
Observe first that if we define the vectors 
\[f_i(t) = e_i + \sum_{k,m}(t-1)\mathbf{g}_2(e_i,h^k_m)h^k_m,\]
then for all $i,j,m$
\[(t\mathbf{g}_1+ (1-t)\mathbf{g}_2)(f_i(t),h^j_m)  = 0\]
We now observe that there must be $\eta^m_{ik}$ such that $d_i = e_i + \sum_{km}\eta^m_{ik}h^k_m$. From the fact that 
$\mathbf{g}_2(d_i,h^k_m) = 0$, 
we deduce that $\eta^m_{ik} = - \mathbf{g}_2(e_i,h^k_m)$.
This shows that $f_i(t) = te_i + (1-t)d_i$ for all $i$, or in other words $f_i(t) = e_i + (t-1)\sum_{km}\eta^m_{ik}h^k_m$, from which it is easy to deduce that  $f_1(t),...,f_n(t),f^1_1,...,f^s_n$ is a polysymplectic basis for all $t$. This polysymplectic basis shows that $t\mathbf{g}_1+ (1-t)\mathbf{g}_2$ is block-compatible with the polysymplectic structure.
\ep
\begin{lem}
\label{existblockcomp}
Let $(M,\omega_1,...,\omega_s)$ be a (non-degenerate) polysymplectic manifold. There exists a Riemannian metric on $M$ block-compatible point by point with the polysymplectic structure.
\end{lem}
{\em Proof}\\
Pick a covering of $M$ by polysymplectic coordinate sets $\mathcal{U}_\alpha$, and a partition of unity $\{f_\alpha\}$ subordinated to the covering.\\
Observe first that if $\mathbf{g}_1$ and $\mathbf{g}_2$ are two Riemannian metrics on $M$ such that for all points $p\in M$ and for any polysymplectic basis $e_1,...,e_n,f^1_1,...,f^s_n$ of $T_pM$, for $j\not = k$, 
$\mathbf{g}_1(f^k_i,f^j_m) = 0 = \mathbf{g}_2(f^k_i,f^j_m) = 0,$
then also $t\mathbf{g}_1+ (1-t)\mathbf{g}_2$ has this property. Therefore,  by using the polysymplectic coordinates on the sets 
$\mathcal{U}_\alpha$, and the partition of unity to sum, we can easily define a Riemannian metric $\mathbf{g}$ on all of $M$ which has the property above at all points $p\in M$.
Define now a family $\mathbf{g}_\alpha$ of block-compatible metrics on any fixed open set  $\mathcal{U}_\alpha$, with the property that $\mathbf{g}_\alpha$ coincides with the fixed $\mathbf{g}$ on the span of $f^1_1,...,f^s_n$ for some, and therefore any, polysymplectic basis. Using the partition of unity, and the previous lemma, we see that we can sum all these metrics to provide a globally defined block-compatible Riemannian metric.
\ep

\begin{lem}
\label{datacompatib}
Let $(M,\omega_1,...,\omega_s)$ be a (non-degenerate) polysymplectic manifold. There is then a one to one correspondence between the following data:\\
1) A Riemannian metric on $M$, compatible with the polysymplectic structure.\\
2) A positive definite non degenerate symmetric bilinear form $\mathbf{g}^1$ on $\bigcap_{j>1}\omega_j^\perp$, plus a constant rank distribution of subspaces $W$ of $TM$, such that at each point $p\in M$ and for some polysymplectic basis $e_1,...,e_n,f^1_1,...,f^s_n$ of $T_pM$, $\mathbf{g}^1|_{T_pM}$ is supported on the span of $f^1_1,...,f^1_n$,  and $W_p = <e_1,...,e_n>$.\\
In the direction from $1)$ to $2)$ the correspondence sends a metric $\mathbf{g}$ to the bilinear form $\mathbf{g}^1$ and the subspace $W$ defined for any $p$ and any polysymplectic basis  $e_1,...,e_n,f^1_1,...,f^s_n$ of $T_pM$ as $\mathbf{g}^1|_{T_pM} = \mathbf{g}|_{<f^1_1,...,f^1_n>}$ and $W_p = <f^1_1,...,f^s_n>^{\perp_\mathbf{g}}$ respectively
\end{lem}
\textit{Proof}\\
In the direction from $1)$ to $2)$, to check that the correspondence is well defined it is enough to observe that $W_p = <e_1,...,e_n>$ for any orthonormal polysymplectic basis  $e_1,...,e_n,f^1_1,...,f^s_n$.
In the direction from $2)$ to $1)$, to define $\mathbf{g}|_{T_pM}$ choose any polysymplectic basis $e_1,...,e_n,f^1_1,...,f^s_n$ such that $W_p = <e_1,...,e_n>$, and 
$f^1_1,...,f^1_n$ is $\mathbf{g}^1|_{T_pM}$-orthonormal. Then declare any such basis to be $\mathbf{g}$-orthonormal. To check that this definition is correct, suppose given any other polysymplectic basis with the same property. Then it is immediate to check, using Lemma \ref{injectivityssymp}, and the observation that if a matrix is orthogonal also the transpose of its inverse is so (and actually coincides with it), that the transition matrix from one basis to the other is orthogonal, and therefore $\mathbf{g}$ is well defined. By construction, the metric $\mathbf{g}$ is Riemannian, and compatible with the polysymplectic structure point by point. The verification that the metric defined varies smoothly as $p$ varies in $M$ is straightforward, and left to the reader.
Both the correspondences thus defined are one to one and onto, as they are one the inverse of the other.
\ep

\textit{Proof of the theorem}\\
Pick any globally defined block-compatible Riemannian metric $\mathbf{g}_0$ on $M$, which exists from Lemma ~\ref{existblockcomp}. At any given point $p\in M$, pick any polysymplectic basis $e_1,...,e_n,f^1_1,...,f^s_n$, and consider the bilinear form $\mathbf{g}^1|_{T_pM} = \mathbf{g}|_{<f^1_1,...,f^1_n>}$ and the subspace $W_p = <f^1_1,...,f^s_n>^{\perp_\mathbf{g}}$. The bilinear form  $\mathbf{g}^1$ and the distribution of subspaces $W$ thus defined determine uniquely a Riemannian metric compatible with the polysymplectic structure, in view of Lemma \ref{datacompatib}.\\
To see that the space of compatible metrics is contractible, pick any metric $\mathbf{g}_0$ in it. Using Lemma \ref{datacompatib}, it is easy to see that there is a canonical way to interpolate between 
$\mathbf{g}_0$ and any other metric $\mathbf{g}$ compatible with the polysymplectic structure, and that this interpolation procedure provides a retraction of the space of compatible metrics to its point $\mathbf{g}_0$.
\ep

\section{Lagrangians and the Legendre transform}
Let's assume that we are in the following setting: there are a target space, $\mathbf{R}^{n}$ (with coordinates $x_{1},...,x_{n}$), and a first order Lagrangian depending on $s$ variables, i.e.
\[\mathcal{L}~=~F\left(x_{i},\frac{\partial{x_{i}}}{\partial z_{j}}\right)\]
Any such Lagrangian can be seen naturally as a map
\[\mathcal{L}~:~\st(\mathbf{R}^{n})~\rightarrow~\mathbf{R}\]
Note that this map is what is classically defined to be a Lagrangian in the case $s~=~1$ (the definition of $\st(M)$ for a smooth manifold $M$ is given in the first section). We formalize these considerations with the following (classical)
\begin{dfn}
Let $M$ be a smooth manifold. A {\em (first order) Lagrangian} on $M$ is a function $\mathcal{L}$ from $\st(M)$ to $\rnum$.
\end{dfn} 
\begin{dfn}
\label{local-form-omegas}
A Lagrangian 
$\mathcal{L}~:~\st(\mathbf{R}^{n})~\rightarrow~\mathbf{R}$
is {\em nondegenerate} at a point $p$ if the forms
\[d\left(\sum_{i=1}^{n}\frac{\partial \mathcal{L}}{\partial y^{j}_{i}}dx_{i}\right),~ ~j~=~1,...,s\]
induce a polysymplectic structure on $\st(\mathbf{R}^{n})_p$. We use the notation
\[\omega_{\mathcal{L}}^{j}~=~\sum_{i=1}^{n}\frac{\partial \mathcal{L}}{\partial y^{j}_{i}}dx_{i}\]
\end{dfn}
We used the symbol $\omega^j$ for the forms above to adhere tothe classical notation used in the symplectic case ($s = 1$). One should not confuse these $1$-forms with the forms $\omega_j$, which are usually two-forms part of a polysymplectic structure.\\ 
To extend the definition from $\rnum^n$ to a general $M$ we need the following
\begin{lem}
Let $\Phi~:~\mathbf{R}^{n}~\rightarrow~\mathbf{R}^{n}$ be a diffeomorphism, and let $\mathcal{L}$ be a s-Lagrangian. If we use the notation $\tilde{x}_i = x_i\tilde{\Phi}$,  $\tilde{y}^j_i$ for the coordinates on $\st(\rnum^n)$ corresponding to the $\tilde{x}_i$ on $\rnum^n$, and $\tilde{\omega}^j_\mathcal{L}$ for  $\sum_{i=1}^{n}\frac{\partial \mathcal{L}}{\partial \tilde{y}^{j}_{i}}d\tilde{x}_{i}$, we have that $\tilde{\omega}^j_\mathcal{L} = \omega_{\mathcal{L}}^{j}$
\end{lem}
{\it Proof}\\
Because  $\left(\frac{\partial}{\partial \tilde{x}_i}\right)_j = \sum_l \frac{\partial x_l}{
\partial \tilde{x}_i}\left(\frac{\partial}{\partial x_l}\right)_j$, we obtain that $\delta_{mi}\delta_{kj} = \tilde{y}^j_i\left(\left(\frac{\partial}{\partial \tilde{x}_m}\right)_k\right) = \sum_l 
\frac{\partial x_l}{\partial \tilde{x}_m}\tilde{y}^j_i\left(\left(\frac{\partial}{\partial x_l}\right)_k\right) $, from which it follows that $\tilde{y}^j_i = \sum_n \frac{
\partial \tilde{x}_i}{\partial x_n}y^j_n$, $~ y^j_i = \sum_n \frac{\partial x_i}{
\partial \tilde{x}_n}\tilde{y}^j_n$. We then obtain that 
\[\tilde{\omega}^j_\mathcal{L} = \sum_{i=1}^{n}\frac{\partial \mathcal{L}}{\partial \tilde{y}^{j}_{i}}d\tilde{x}_{i} =
\sum_{i=1}^{n}\left(\sum_{m =1}^n \frac{\partial \mathcal{L}}
{\partial x_{m}}
\frac{\partial x_{m}}{\partial \tilde{y}^{j}_{i}} + \sum_{k,m}\frac{\partial \mathcal{L}}
{\partial y^{k}_{m}}
\frac{\partial y^{k}_{m}}{\partial \tilde{y}^{j}_{i}}\right)\sum_{l = 1}^n\frac{\partial \tilde{x}_{i}}{\partial x_l}dx_l =\]
\[ \sum_{i=1}^{n}\left(\sum_{m =1}^n\frac{\partial \mathcal{L}}
{\partial y^{j}_{m}}
\frac{\partial x_m}{
\partial \tilde{x}_i}\right)\sum_{l = 1}^n\frac{\partial \tilde{x}_{i}}{\partial x_l}dx_l = \sum_{m =1}^n\frac{\partial \mathcal{L}}
{\partial y^{j}_{m}}dx_m = \omega^j_\mathcal{L}\]
\ep

\begin{exe}
Let $\epsilon(i)\in\{-1,1\}$ for $i=1,...,s$, and let
\[\mathcal{L}~=~\frac{1}{2}\sum_{i,j}\epsilon(j)(y^{j}_{i})^{2}\]
Then
\[\omega_{\mathcal{L}}^{j}~=~\epsilon(j)\sum_{i}y^{j}_{i}dx_{i},~ ~
d\omega_{\mathcal{L}}^{j}~=~\epsilon(j)\sum_{i}dy^{j}_{i}\wedge dx_{i}\]
At this point it is clear that we obtained a polysymplectic structure, and therefore the $s$-Lagrangian is non degenerate at every point.
\end{exe}
\begin{dfn}
Let $M$ be a smooth manifold of dimension $n$. Then an {\em $s$-Lagrangian} (or simply a Lagrangian) on $M$ is a smooth map
\[\mathcal{L}~:~\st(M)~\rightarrow~\mathbf{R}\]
\end{dfn}
\begin{teo}[Definition of $\omega^j_\mathcal{L}$ and $H_\mathcal{L}$]
Given a smooth manifold $M$ and a Lagrangian
$\mathcal{L}~:~\st(M)~\rightarrow~\mathbf{R}$:\\
1) There is a unique $\mathbf{Diff}(M)$-equivariant choice of $s$ forms
\[\omega^j_\mathcal{L}~\in~\omega_{1}\left(\st(M)\right),~ ~j~=~1,...,s\]
which in local coordinates are  expressed as in Definition ~\ref{local-form-omegas}.\\
2) There is a unique function $H_\mathcal{L}$ on $\st(M)$ which is expressed locally as $\sum_{m,j}y^j_m\frac{\partial \mathcal{L}}{\partial y^j_m}~-~\mathcal{L}$, for any choice of local coordinates.
\end{teo}
{\it Proof}\\
1) The previous Lemma shows that there is a canonical choice for the $\omega^j_{\mathcal{L}}$. Namely, for $p\in M$, choose a system of coordinates $x_1,...,x_n$ in a neighborhood $\mathcal{U}\subset M$ of $p$, and let $x_1,...,x_n,y^1_1,..,y^s_n$ be the corresponding coordinates on $\st(\mathcal{U})\subset\st(M)$. We then take 
\[\omega_{\mathcal{L}}^{j}~=~\sum_{i=1}^{n}\frac{\partial \mathcal{L}}{\partial y^{j}_{i}}dx_{i}\] 
The previous Lemma guarantees that this definition does not depend on the choice of the coordinates $x_1,...,x_n$. To see that the forms are $\mathbf{Diff}(M)$-equivariant, let $\Phi~:M\to M$ be a diffeomorphism. On $\Phi^{-1}(\mathcal{U})$, take coordinates $\tilde{x}_1,...,\tilde{x}_n$ with $\tilde{x}_i = x_i\Phi$, and let $\tilde{x}_1,...,\tilde{x}_n,\tilde{y}^1_1,..,\tilde{y}^s_n$ be the corresponding coordinates on $\st(\Phi^{-1}(\mathcal{U}))\subset\st(M)$. It is immediate to check that $\tilde{y}^j_i = y^j_i\Psi$, where $\Psi : \st(M)\to\st(M)$ is the canonical extension of $\Phi$. Let $\tilde{\mathcal{L}} = \mathcal{L}\Psi$. With these coordinates, the map on $\st(\rnum^n)$ associated to $\Psi$ is just the identity. Therefore, we can conclude immediately that $\Psi^*\omega^j_\mathcal{L} = \sum_{i=1}^{n}\frac{\partial \tilde{\mathcal{L}}}{\partial \tilde{y}^{j}_{i}}d\tilde{x}_{i}$ and from the uniqueness property we conclude that $\Psi^*\omega^j_\mathcal{L} = \omega^j_{\mathcal{L}\Psi}$.\\
2) Using the computations make in the previous Lemma, we know that if $x_1,...,x_n$ and $\tilde{x}_1,...,\tilde{x}_n$ are two systems of coordinates on (a neighborhood of a point in) $M$, and the $x_i,y^j_i$, $\tilde{x}_i,\tilde{y}^j_i$ are the corresponding coordinates on $\st(M)$, we have $\tilde{y}^j_i = \sum_n \frac{
\partial \tilde{x}_i}{\partial x_n}y^j_n$, $~ y^j_i = \sum_n \frac{\partial x_i}{
\partial \tilde{x}_n}\tilde{y}^j_n$. If $\tilde{H}_\mathcal{L}$ is the expression obtained from the $\tilde{x}_i,\tilde{y}^j_i$ coordinates,
\[\tilde{H}_\mathcal{L} + \mathcal{L} = \sum_{m,j}\tilde{y}^j_m\frac{\partial \mathcal{L}}{\partial \tilde{y}^j_m}  = \sum_{m,j}\sum_{n,l}\frac{
\partial \tilde{x}_m}{\partial x_n}y^j_n\frac{\partial \mathcal{L}}{\partial y^j_l}\frac{\partial x_l}{
\partial \tilde{x}_m} = \sum_{j,l}
y^j_l\frac{\partial \mathcal{L}}{\partial y^j_l} = H_\mathcal{L} + \mathcal{L}\]
\ep
\begin{dfn}
A smooth lagrangian $\mathcal{L}~:~\st(M)\to \rnum$ is {\em nondegenerate at the point $p\in \st(M)$} if the forms $d\omega^1_\mathcal{L},...,d\omega^s_\mathcal{L}$ induce a polysymplectic structure on $\st(M)$ at $p$. It is  {\em nondegenerate} if it is nondegenerate at all points.
\end{dfn}
Note that it is very easy to prove that if a lagrangian $\mathcal{L}~:~\st(M)\to \rnum$ is nondegenerate, the forms $d\omega^1_\mathcal{L},...,d\omega^s_\mathcal{L}$ induce a polysymplectic structure on $\st(M)$.
\begin{dfn}
Given a smooth manifold $M$ and a $\mathbf{C}^{1}$ map $f~:~\mathbf{R}^{s}~\rightarrow~M$, there is a naturally induced map
\[\tilde{f}~:~\mathbf{R}^{s}~\rightarrow~\st(M)\]
obtained by taking
$\tilde{f}(z)~=~\left(f(z),df_z\left(\frac{\partial}{\partial z_{1}}\right),...,
df_z\left(\frac{\partial}{\partial z_{s}}\right)\right)$, where $df_z$ is the differential of $f$ at the point of coordinates $(z_1,...,z_s)$ of $\rnum^s$, and the $\frac{\partial}{\partial z_{i}}$ are the standard basis for $T_z\rnum^s$
\end{dfn}
\begin{lem}
Given a smooth manifold $M$, a lagrangian $\mathcal{L}$ on it and a $\mathbf{C}^{2}$ map $f~:~\mathbf{R}^{s}~\rightarrow~M$, we have that in local coordinates $x_i,y^j_i$ on $\st(M)$,
\[d\tilde{f}_z(\frac{\partial}{\partial z_{k}})\contr d\omega_{\mathcal{L}}^{j}~=~\sum_i\frac{\partial}{\partial z_k}\left(\frac{\partial ~\mathcal{L}}{\partial y^j_i}(\tilde{f}(z))\right)dx_i - \sum_{l,m,i}\frac{\partial^2 ~\mathcal{L}}{\partial y^l_m\partial y^j_i}\frac{\partial f_i}{\partial z_{k}}dy^l_m\]
\[  - \sum_{m,i}\frac{\partial^2 ~\mathcal{L}}{\partial x_m\partial y^j_i}\frac{\partial f_i}{\partial z_{k}}dx_m\]
\end{lem}
\textit{Proof} This is just a direct computation.\\
$d\tilde{f}_z(\frac{\partial}{\partial z_{k}}) = \sum_{l,m}\frac{\partial ~y^l_m(\tilde{f}(z))}{\partial z_{k}}\frac{\partial}{\partial y^l_m} + \sum_i\frac{\partial f_i}{\partial z_{k}}\frac{\partial}{\partial x_i} = \sum_{l,m}\frac{\partial^2 ~f_m(z))}{\partial z_{k}\partial z_l}\frac{\partial}{\partial y^l_m} + \sum_i\frac{\partial f_i}{\partial z_{k}}\frac{\partial}{\partial x_i}$.\\
$d\tilde{f}_z(\frac{\partial}{\partial z_{k}})\contr d\omega_{\mathcal{L}}^{j} = \sum_{l,m,i}\frac{\partial^2 ~\mathcal{L}}{\partial y^l_m\partial y^j_i}\left(\frac{\partial^2 ~f_m}{\partial z_{k}\partial z_l}dx_i - \frac{\partial f_i}{\partial z_{k}}dy^l_m\right) +$\\
$\sum_{m,i}\frac{\partial^2 ~\mathcal{L}}{\partial x_m\partial y^j_i}\left(\frac{\partial f_m}{\partial z_{k}}dx_i - \frac{\partial f_i}{\partial z_{k}}dx_m\right) = \sum_i\frac{\partial}{\partial z_k}\left(\frac{\partial ~\mathcal{L}}{\partial y^j_i}(\tilde{f}(z))\right)dx_i -$\\ 
$\sum_{l,m,i}\frac{\partial^2 ~\mathcal{L}}{\partial y^l_m\partial y^j_i}\frac{\partial f_i}{\partial z_{k}}dy^l_m - \sum_{m,i}\frac{\partial^2 ~\mathcal{L}}{\partial x_m\partial y^j_i}\frac{\partial f_i}{\partial z_{k}}dx_m$.
\ep
\begin{teo}
\label{ele-slag}
Let $M$ be a smooth manifold , let $\mathcal{L}$ be a lagrangian  on it and  $f~:~\mathbf{R}^{s}~\rightarrow~M$ a $\mathbf{C}^{2}$ map. If we indicate with $\tilde{f}$ the natural lifting of $f$ to $\st(M)$, and with $z_1,...,z_s$ the canonical coordinates on $\rnum^s$, and $z^j$ on $\mathbf{R}^{s}$, the following are equivalent:\\
1) For any choice of  local coordinates $x_1,...,x_n$ on $M$ around a point in the image of $f$, the map $f$ satisfies the Euler-Lagrange equations with respect to $\mathcal{L}$, 
\[ \frac{\partial \mathcal{L}}{\partial x_i}(\tilde{f}(z)) - \sum_{j=1}^s \frac{\partial}{\partial z_{j}}\left(\frac{\partial \mathcal{L}}{\partial y^j_i}(\tilde{f}(z))\right) = 0~ ~ ~\forall i\in\{1,...,n\}\]
2) \[\sum_{j=1}^s\left(d\tilde{f}\left(\frac{\partial}{\partial z_j}\right)\contr d\omega^j_\mathcal{L}\right)|_{\tilde{f}(\mathbf{z})}~=~-dH_\mathcal{L}|_{\tilde{f}(\mathbf{z})}\]
\end{teo}
\textit{Proof} The proof is just a direct computation. 
First of all, observe that\\ 
$y^j_m(\tilde{f}) = \left(\frac{\partial}{\partial x_m}\right)^*\left(df\left(\frac{\partial}{\partial z_j}\right)\right) = \left(\frac{\partial}{\partial x_m}\right)^*\left(\sum_l \frac{\partial f_l}{\partial z_j}\frac{\partial}{\partial x_l}\right) = \frac{\partial f_m}{\partial z_j}$.\\ 
We then have that\\
${dH_\mathcal{L}}|_{\tilde{f}(\mathbf{z})} = \sum_i\left(\sum_{m,j} \frac{\partial f_m}{\partial z_j}\frac{\partial^2 \mathcal{L}}{\partial x_i\partial y^j_m} - \frac{\partial \mathcal{L}}{\partial x_i}\right)dx_i + \sum_{j,m,h,l}\frac{\partial^2 \mathcal{L}}{\partial y^h_l\partial y^j_m}\frac{\partial f_m}{\partial z_j}dy^h_l$, and\\
$\sum_{j=1}^s\left(d\tilde{f}\left(\frac{\partial}{\partial z_j}\right)\contr d\omega^j_\mathcal{L}\right)|_{\tilde{f}(\mathbf{z})} = \sum_{i,j}\left(\frac{\partial}{\partial z_j}\left(\frac{\partial \mathcal{L}}{\partial y^j_i}(\tilde{f})\right) - \sum_{m}\frac{\partial^2 \mathcal{L}}{\partial x_i\partial y^j_m}\frac{\partial f_m}{\partial z_j}\right)dx_i -$\\ 
$\sum_{h,m,i,j}\frac{\partial^2 \mathcal{L}}{\partial y^h_m\partial y^j_i}\frac{\partial f_i}{\partial z_j}dy^h_m$. Therefore,\\
$\sum_{j=1}^s\left(d\tilde{f}\left(\frac{\partial}{\partial z_j}\right)\contr d\omega^j_\mathcal{L}\right)|_{\tilde{f}(\mathbf{z})} + {dH_\mathcal{L}}|_{\tilde{f}(\mathbf{z})} =
\sum_i\left(- \frac{\partial \mathcal{L}}{\partial x_i} + \sum_{j}\frac{\partial}{\partial z_j}\left(\frac{\partial \mathcal{L}}{\partial y^j_i}(\tilde{f})\right)\right)dx_i$ and the thesis follows.
\ep
\begin{rmk}
Given a smooth manifold $M$ and a Lagrangian
$\mathcal{L}~:~\st(M)~\rightarrow~\rnum$, 
the forms $\omega_{\mathcal{L}}^{j}$ on $\st(M)$ are semi-basic, 
i.e. for any $v\in\st(M)$ and any $j\in\{1,...,s\}$ 
there is a (necessarily unique) covector $\tau_{\mathcal{L}}^{j}(v)\in\mathbf{T}^*M_{\pi(v)}$ 
such that 
\[\pi_j^*\left(\tau_{\mathcal{L}}^{j}(v)\right)~=~(\omega_{\mathcal{L}}^{j})_v\]
\end{rmk}

\begin{dfn}
Given a smooth manifold $M$ and a Lagrangian
$\mathcal{L}~:~\st(M)~\rightarrow~\rnum$, the {\em Legendre transformation} from  $\st(M)$ to $\st^*(M)$ (relative to $\mathcal{L}$) is defined for $v\in\st(M)$ as
$\tau_{\mathcal{L}}(v)~:=~\oplus_{j=1}^s \tau_{\mathcal{L}}^{j}(v)$.
\end{dfn}
\begin{teo}
\label{legendre-ssymp}
Given a smooth manifold $M$ and a Lagrangian
$\mathcal{L}~:~\st(M)~\rightarrow~\rnum$, the Legendre transformation induces a morphism (also called Legendre transformation)
\[(\st(M),d\omega_{\mathcal{L}}^{1},...,d\omega_{\mathcal{L}}^{s})\to (\st^*(M),\omega_1,...,\omega_s)\]
(where the $\omega_1,...,\omega_s$ are the forms introduced in Example ~\ref{execotangent}), in the sense that $\tau_{\mathcal{L}}^*\omega_j = d\omega^j_\mathcal{L}$. If the lagrangian is nondegenerate, the Legendre transformation is a local diffeomorphism, and a morphism of polysymplectic manifolds.
\end{teo}
\textit{Proof}
It is clearly enough to prove the theorem locally (on $M$). We may therefore assume without loss of generality that $M=\rnum^d$. We take coordinates $x_i$ on $\rnum^d$, $x_i,y^j_i$ on $\st(\rnum^d)$ and $x_i,q^j_i$ on $\st^*(\rnum^d)$. In terms of these coordinates, $\omega_{\mathcal{L}}^{j}~=~\sum_{i}\frac{\partial \mathcal{L}}{\partial y^{j}_{i}}dx_{i}$ and $\omega_j~=~\sum_i dq^j_i\wedge dx_i$. We have
\[d\omega_{\mathcal{L}}^{j}~=~\sum_{k,l}\sum_{i}\frac{\partial^2 \mathcal{L}}{\partial y^{k}_{l}\partial y^{j}_{i}}dy^{k}_{l}\wedge dx_{i}~+~
\sum_{m}\sum_{i}\frac{\partial^2 \mathcal{L}}{\partial x_m\partial y^{j}_{i}}dx_m\wedge dx_{i}\]
We have therefore that 
\[\tau_{\mathcal{L}}(v)~=~(\sum_{i}\frac{\partial \mathcal{L}}{\partial y^{1}_{i}}(v)dx_{i},...,\sum_{i}\frac{\partial \mathcal{L}}{\partial y^{s}_{i}}(v)dx_{i})\]
or in other words $q^j_i\left(\tau_{\mathcal{L}}(v)\right) = \frac{\partial \mathcal{L}}{\partial y^{j}_{i}}(v)$, $x_i\left(\tau_{\mathcal{L}}(v)\right) = x_i(v)$. By definition, $\tau_{\mathcal{L}}^*(\omega_j)=\sum_i d(\frac{\partial \mathcal{L}}{\partial y^{j}_{i}})\wedge dx_i$, from which it apparent that $\tau_{\mathcal{L}}^*(\omega_j)~=~d\omega_{\mathcal{L}}^{j}$.\\
If $\mathcal{L}$ is non degenerate, to prove that $\tau_{\mathcal{L}}$ is a local diffeomorphism we will show that $d\tau_{\mathcal{L}}$ is injective (and therefore an isomorphism, by dimension count) at every point.
\[d\tau_{\mathcal{L}}(v)\left(\sum_{k,l}\alpha_k^l\frac{\partial}{\partial y^k_l}~+~\sum_m\beta^m\frac{\partial}{\partial x_m}\right)~=~0\]
implies $\sum_{k,l}\alpha_k^l\frac{\partial^2 \mathcal{L}}{\partial y^{k}_{l}\partial y^{j}_{i}}~=~0~\forall i,j$ and $\beta_m~=~0 \forall m$. The condition on the $\alpha_k^l$ is easily seen to be equivalent to
$\forall j ~ ~d\omega_{\mathcal{L}}^{j}\contr \left(\sum_{k,l}\alpha_k^l\frac{\partial}{\partial y^k_l}\right)~=~0$. To conclude, we observe that the fact that any vector which contracts to zero with all the structure forms must be zero is a basic property of polysymplectic vector spaces, proved in Lemma ~\ref{injectivityssymp}.
\ep
\begin{dfn}
1) An s-Lagrangian $\mathcal{L}$ on a manifold $M$ is said to be {\em strongly non degenerate} if the associated Legendre transformation is a smooth diffeomorphism.\\
2) Given a smooth manifold $M$ and a strongly non degenerate Lagrangian
$\mathcal{L}$, the function $\mathcal{H}(\mathcal{L})$ on $\st^*(M)$, called the {\em Hamiltonian} associated to $\mathcal{L}$, is defined as  $\mathcal{H}(\mathcal{L}) = H_\mathcal{L}\tau_\mathcal{L}^{-1}$. We sometimes write $\mathcal{H}(\mathcal{L}) = \sum_{i,j}y^j_iq^j_i - \mathcal{L}$, meaning that all the functions defined on $\st(M)$ are considered as functions on $\st^*(M)$ via the Legendre transformation.\\
\end{dfn}
Note that the condition of strong non degeneracy can be somewhat weakened, if we restrict the domain of the Legendre transformation to some open domain in $\st(M)$ which is not necessarily of the form $\st(\mathcal{U})$, with $\mathcal{U}\subset M$. We will not pursue this in full generality for the moment, and we will content ourselves with the study of the situation described in the following theorem. The brackets $\{~,~\}_j$ on $\st^*(M)$ were defined in Theorem ~\ref{canspoisson}.
\begin{teo}
\label{gener-ham}
Given a smooth manifold $M$ and a Lagrangian
$\mathcal{L} : \st(M)\to\rnum$, assume that the Legendre transformation associated to $\mathcal{L}$ is invertible on the open set $\mathcal{V}\subset \st(M)$. Let $\mathcal{H} = \mathcal{H}(\mathcal{L}) = H_\mathcal{L}\tau_\mathcal{L}^{-1}$. Assume moreover that we have chosen coordinates $x_1,...,x_n$ on $M$ (possibly restricting from $M$ to a proper open subset), and that the $x_1,...,x_n,q^1_1,...,q^s_n$ are the associated coordinates on $\st^*(M)$. Then, if $f~:~\mathcal{W}\subset\mathbf{R}^{s}~\rightarrow~M$ is a $\mathbf{C}^{2}$ map with $\tilde{f}\left(\mathcal{W}\right)\subset \mathcal{V}$, the following are equivalent:\\
1) $f$ satisfies the Euler-Lagrange equations with respect to $\mathcal{L}$.\\
2) If $q^j_i(z):=~q^J_i(\tau_\mathcal{L}(\tilde{f}(z)))$, $x_i~:=~x_i(\tau_\mathcal{L}(\tilde{f}(z)))$, then
\[\forall i~ ~\sum_j\frac{\partial}{\partial z_j}(q^j_i(z)) = \frac{1}{s}\sum_j\{ q^j_i,\mathcal{H}\}_j~\left(= -\frac
{\partial \mathcal{H}}
{\partial x_i}\right),~ ~ ~\forall k,i~ ~\frac{\partial x_i(z)}{\partial z_k} = \left\{x_i(z),\mathcal{H}\right\}_k\]
\end{teo}
\textit{Proof}\\
From Theorem ~\ref{ele-slag} we know that the function $f$ satisfies the Euler-Lagrange equations if and only if 
$\sum_{j=1}^s\left(d\tilde{f}\left(\frac{\partial}{\partial z_j}\right)\contr d\omega^j_\mathcal{L}\right)|_{\tilde{f}(\mathbf{z})}~=~-dH_\mathcal{L}|_{\tilde{f}(\mathbf{z})}$.
If we indicate with $\tilde{\tilde{f}}$ the map $\tau_\mathcal{L}\tilde{f}$, and we push forward the equation above using $\tau_\mathcal{L}$, it becomes, using the Leibnitz rule and Theorem ~\ref{legendre-ssymp},
\[\sum_{j=1}^s\left(d\tilde{\tilde{f}}\left(\frac{\partial}{\partial z_j}\right)\contr \omega_j\right)|_{\tilde{\tilde{f}}(\mathbf{z})}~=~-d\mathcal{H}|_{\tilde{\tilde{f}}(\mathbf{z})}\]
We have also that\\
$\sum_{j=1}^s\left(d\tilde{\tilde{f}}\left(\frac{\partial}{\partial z_j}\right)\contr \omega_j\right)|_{\tilde{\tilde{f}}(\mathbf{z})} = \sum_{j=1}^s\left(\left(\sum_i\frac{\partial x_i(\mathbf{z})}{\partial z_j}dx_i + \sum_{h,m}\frac{\partial q^h_m(\mathbf{z})}{\partial z_j}dq^h_m\right)\contr \omega_j\right)$\\
$= \sum_{j=1}^s\sum_{i=1}^n\left(\frac{\partial x_i(\mathbf{z})}{\partial z_j}dq^j_i - \frac{\partial q^j_i(\mathbf{z})}{\partial z_j}dx_i\right)$, and hence\\
$\sum_{j=1}^s\left(d\tilde{\tilde{f}}\left(\frac{\partial}{\partial z_j}\right)\contr \omega_j\right)|_{\tilde{\tilde{f}}(\mathbf{z})} + d\mathcal{H}|_{\tilde{\tilde{f}}(\mathbf{z})} =
\sum_{j=1}^s\sum_{i=1}^n\left(\frac{\partial x_i(\mathbf{z})}{\partial z_j} + 
\frac
{\partial \mathcal{H}}
{\partial q^j_i}\right)dq^j_i +$\\
$\sum_{i=1}^n\left(- \sum_j\frac{\partial q^j_i(\mathbf{z})}{\partial z_j} + \frac
{\partial \mathcal{H}}
{\partial x_i}\right)dx_i = $\\
$\sum_{i,j}\left(\frac{\partial x_i(\mathbf{z})}{\partial z_j} + 
\{x_i,\mathcal{H}\}_j\right)dq^j_i + \sum_{i,j}\left(- \frac{\partial q^j_i(\mathbf{z})}{\partial z_j} - \frac{1}{s}\{q^j_i,\mathcal{H}\}_j\right)dx_i$
\ep

\section{Almost s-\ka manifolds and Special Lagrangian Fibrations}
\label{sec:almostska}
\begin{dfn}
\label{dfn-almostska}
A  smooth polysymplectic manifold together with a Riemannian metric 
$(M,\omega_{1},...,\omega_{s},\mathbf{g})$ 
is {\em ~almost s-\ka~} if the metric is compatible with the 
polysymplectic structure
\end{dfn}
\begin{exe}[With assigned topology]
\label{cotan-almostska}
Let $(M,\mathbf{g})$ be a Riemannian manifold. Then $\st^{*}(M)$, 
together with its canonical polysymplectic structure and the induced 
metric, is almost $s$-\ka.
\end{exe}
\textit{Proof}
We consider the case $s=1$. The general case can be done similarly.\\
As the question is local, we may assume $M=\rnum^n$. We assume also that we have fixed coordinates $x_1,...,x_n$ on $M$ around $0$, so that the metric $\me$ is given at any point by $\me_{ij} = \me(\frac{\partial}{\partial x_i},\frac{\partial}{\partial x_j})$. The vector fields $\frac{\partial}{\partial x_1},...,\frac{\partial}{\partial x_n}$ induce coordinates on $T^*M$, which we will indicate with $y_1,...,y_n$. Therefore, the point on $T^*M$ corresponding to coordinates $(x,y)$ will be $(p_x,\sum_{i=1}^ny_idx_i)$. As a general rule, we will use $1,...,n$ as indices for the $x$ coordinates, and $\bar{1},...,\bar{n}$ as indices for the $y$ coordinates.
A change in coordinates on $M$ from the $x$ to the new coordinates $\tilde{x}$ induces automatically a change from the $(x,y)$ coordinates to the $\tilde{x},\tilde{y}$ coordinates with $\tilde{y}_i = \sum_jy_j\frac{\partial x_j}{\partial \tilde{x}_i}$, as $\sum_{i=1}^ny_idx_i = \sum_{i,j}y_i\frac{\partial x_i}{\partial \tilde{x}_j}d\tilde{x}_j$. The tangent frame $\frac{\partial}{\partial x_1},...,\frac{\partial}{\partial x_n},\frac{\partial}{\partial y_1},...,\frac{\partial}{\partial y_n}$  at the point corresponding to $(x,y)$ is therefore transformed by any change of coordinates of the above form as follows:
\[\frac{\partial}{\partial x_k} = \sum_i\frac{\partial \tilde{x}_i}{\partial x_k}\frac{\partial}{\partial \tilde{x}_i} + \sum_{i}\frac{\partial \tilde{y}_i}{\partial x_k}\frac{\partial}{\partial \tilde{y}_i},~ ~ ~\frac{\partial}{\partial y_k} = \sum_i\frac{\partial x_k}{\partial\tilde{x}_i}\frac{\partial}{\partial \tilde{y}_i}\]
However, $\frac{\partial \tilde{y}_i}{\partial x_k} = \frac{\partial}{\partial x_k}\left(\sum_jy_j\frac{\partial x_j}{\partial \tilde{x}_i}\right) = \sum_jy_j\frac{\partial}{\partial x_k}\left(\frac{\partial x_j}{\partial \tilde{x}_i}\right) = \sum_{j,l}y_j\frac{\partial \tilde{x}_l}{\partial x_k}\frac{\partial^2 x_j}{\partial \tilde{x}_i\tilde{x}_l}$
And therefore the first $n$ vectors are transformed as
\[\frac{\partial}{\partial x_k} = \sum_i\frac{\partial \tilde{x}_i}{\partial x_k}\left(\frac{\partial}{\partial \tilde{x}_i} + \sum_{j,l}y_j\frac{\partial^2 x_j}{\partial \tilde{x}_l\partial \tilde{x}_i}\frac{\partial}{\partial \tilde{y}_l}\right)\]
We define a new frame at the point corresponding to coordinates $(x,y)$ as
\[v_k = \frac{\partial}{\partial x_k} + \sum_{l,m}y_l\Gamma^l_{k m}\frac{\partial}{\partial y_m},~ ~ ~ ~\frac{\partial}{\partial y_k};~ ~ ~ k=1,...,n\]
where the $\Gamma^l_{k m}$ are Christoffel's symbols relative to the metric $\me$. We then have that, using the transformation laws for Christoffel's symbols, the frame above transforms under a coordinate change of the form described above as
\[\tilde{v}_k = \sum_i\frac{\partial \tilde{x}_i}{\partial x_k}v_i,~ ~ ~\frac{\partial}{\partial y_k} = \sum_i\frac{\partial x_k}{\partial\tilde{x}_i}\frac{\partial}{\partial \tilde{y}_i}\]
From these transformation laws we see that we can define a metric structure $\eme$ on $T^*(M)$ by giving the metric tensor on the above frame as follows:
\[\eme(v_i,v_j) = \me_{ij},~ ~\eme(\frac{\partial}{\partial y_i},\frac{\partial}{\partial y_j}) = \me^{ji},~ ~\eme(v_i,\frac{\partial}{\partial y_j}) = 0\]
From this and the definition of the $v_k$ we finally get, using indices $1,...,n$ for the $x$ variables and $\bar{1},...,\bar{n}$ for the $y$ variables,
\[\eme_{ij} = \me_{ij} + \sum_{lmhr}y_jy_r\Gamma^l_{im}\Gamma^r_{jh}\me^{mh},~ ~\eme_{\bar{i}\bar{j}} = \me^{ji},~ ~\eme_{i\bar{j}} = -\sum_{lm}y_l\Gamma^l_{im}\me^{mj}\]
At this point we see already that if we choose the $x$ variables to be a Riemann normal coordinate system around the origin, the coordinate frame induced by the $(x,y)$ variables will be orthonormal on any point on $T^*(M)$ with vanishing $x$ coordinate, and therefore the symplectic form is compatible with the metric.
\ep

The following easy lemma can be used to build many more explicit examples of almost $s$-\ka manifolds
\begin{lem}
Let $\mathbf{G}$ be a connected real Lie group, let $\omega_1,...,\omega_s$ be a polysymplectic structure on it, and assume that all the $\omega_j$ are left-invariant differential forms. Then there is a left-invariant Riemannian metric $\me$ on $\mathbf{G}$ such that $(\mathbf{G},\omega_1,...,\omega_s,\me)$ is  almost $s$-\ka.
\end{lem}
\textit{Proof}\\
As the forms $\omega_j$ define a polysymplectic structure, in particular they induce a polysymplectic structure on the tangent space $T_gG$ for all $g\in G$. Therefore, there exists a polysymplectic basis $\mathcal{V} = (v_1,...,v_n,w^1_1,...,w^s_n)$ of $T_eG$. Then the unique left invariant Riemannian metric $\me_\mathcal{V}$ for which $\mathcal{V}$ is an orthonormal basis of $T_eG$ is easily seen to be compatible with the polysymplectic structure induced on $T_gG$ by the $\omega_j$, for all $g\in G$. By definition,  $(\mathbf{G},\omega_1,...,\omega_s,\me_\mathcal{V})$ is  almost $s$-\ka.
\ep

We now show that starting with $s$  special lagrangian fibrations which are also Riemannian submersions, we can get an almost $s$-\ka manifold. We actually prove a more general result, which then can be specialized to to the setting just mentioned.

The following conditions on a submersion have been already considered in the literature:
\begin{dfn}
Let $(X,\mathbf{g})$ be a Riemannian manifold, let $B$ be a smooth manifold, and let $f~:~X\to B$ be a smooth submersion.\\
1) We say that $f$ is {\em conformal} if there exists a (necessarily unique) conformal structure on $B$ such that $df$ is a conformal map from $Ker(df)_p^\perp$ to $T_{f(p)}B$ for all $p\in X$.\\
2) We say that $f$ is {\em Riemannian} if there exists a (necessarily unique) Riemannian metric on $B$ such that $df$ is an isometry from $Ker(df)_p^\perp$ to $T_{f(p)}B$ for all $p\in X$.\\
3) We say that $f$ is {\em covariant constant} if it is Riemannian, and $df$ commutes with parallel transport, i.e. if $\gamma(t), t\in [0,1]$ is a path in $X$, $G_X~:T_{\gamma(0)}X\to T_{\gamma(1)}X$ is parallel transport in $X$ along $\gamma$, and $G_B~:T_{f(\gamma(0))}B\to T_{f(\gamma(1))}B$ is parallel transport in $B$ along $f(\gamma)$, then 
$G_B\left(df_{\gamma(0)})(v)\right) = df_{\gamma(1)}\left(G_X(v)\right)$ for all $v\in T_pX$.
\end{dfn}
We are now ready to state the main theorem of this section:
\begin{teo}
\label{specialtoska}
Let $(X_i,\omega_{X_i},\mathbf{g}_{X_i},\Omega_{X_i})$  be  \ka  manifolds of complex dimension $n$, for $i = 1,...,s$. Let $B$ be a smooth manifold, and let $f_i:X_i\to B$  be lagrangian fibrations (with respect to the \ka forms) with connected fibres.  Consider $M = X_1\times_B\cdots\times_B X_s$, with the metric $\mathbf{g}$ induced from $X_1\times\cdots\times X_s$ and with the 2-forms $(\omega_1,...,\omega_s)$, where $\omega_i$ is $\sqrt{s}$ times the pull-back of the \ka form of $X_i$, under the natural projection  $M\to X_i$. We then have that:\\
1) $(M,\omega_1,...,\omega_s)$ is a polysymplectic manifold.\\
2) If all the $f_i$ are conformal with respect to the same conformal structure on $B$, then $\mathbf{g}$ is block-compatible with the polysymplectic structure $\omega_1,...,\omega_s$.\\
3)  If all the $f_i$ are Riemannian with respect to the same metric on $B$, then $\mathbf{g}$ is compatible with the polysymplectic structure $\omega_1,...,\omega_s$. In other words,  $(M,\omega_1,...,\omega_s,\mathbf{g})$ is an almost $s$-\ka manifold.\\
4)  If all the $f_i$ are covariant constant with respect to the same metric on $B$, then  all the $\omega_j$ are covariant constant with respect to the metric $\mathbf{g}$ on $M$.
\end{teo}
\textit{Proof}\\
1) This has already been proven in Example ~\ref{teopolysymp}.\\
2) Given $p\in M$, we will show that there is an orthogonal polysymplectic basis of $T_pM$ (which is actually a bit more than bloc-compatibility). Pick an orthogonal basis $v_1,...,v_n$ of $T_{f(p)}B$, and let $z_1^j,...,z_n^j$ be a set of vectors in $Ker(d(f_j)_{p_j})^\perp$  (were $p = (p_1,...,p_s) \in M \subset X_1\times\cdots\times X_s$), such that $df_j(z_i^j) = v_i$ for all $i,j$. Because the $f_j$ are Conformal, it follows that the $z_i^j$ are orthogonal (for fixed $j$). Define 
\[w_i = \frac{1}{\sqrt{|z_i^1|^2 + \cdots + |z_i^s|^2}}(z_i^1,...,z_i^s) \in T_p(X_1\times\cdots\times X_s)\] 
From their definition, it follows that the $w_i$ lie actually in $T_pM$. Moreover, $(w_l,w_m) = \delta_{lm}$. Define also 
\[w_i^j = |z_i^j|(0,...,Jz_i^j,0,..0)~ ~(j^{th}~ place)\]
in $T_pM$. We are indicating with $J$ the complex structure on the various $X_i$ (or the one on $X_1\times\cdots\times X_s$, which is the same). The fact that the $w_i^j \in T_pM$ follows from the fact that $Jz_i^j\in Ker(d(f_j)_{p_j})$, which is a consequence of the Lagrangian condition. It is now very easy to verify that $w_1,...,w_n,w_1^1,...,w_n^s$ is an orthogonal polysymplectic basis at $p$ with respect to the polysymplectic structure $\omega_1,...,\omega_s$.\\
3) Given $p\in M$, we must show that there is an orthonormal polysymplectic basis of $T_pM$. 
The construction of the previous point will give the desired orthonormal basis, provided that we start with  an orthonormal basis $v_1,...,v_n$ of $T_{f(p)}B$.\\
4) As the forms $\omega_j$ are clearly covariant constant on $X_1\times\cdots\times X_s$ (because they are \ka forms, and hence covariant constant on their respective $X_j$'s), it is enough to observe that if all the $f_j$ are covariant constant, then $M$ is a totally geodesic submanifold of  $X_1\times\cdots\times X_s$. Indeed, parallel transport on $M$ is then just the restriction of parallel transport on $X_1\times\cdots\times X_s$, and hence the $\omega_j$ are constant also on $M$.

\ep

To put the condition of being Riemannian into perspective, we connect it with the semi-flatness condition of ~\cite{SYZ}, or rather with one of its consequences. We start by recalling the following standard
\begin{dfn}
Let $(X,\omega,\mathbf{g},\Omega)$ be a Calabi-Yau  manifold of complex dimension $n$ (where $\omega$ is the \ka form, $\mathbf{g}$ the \ka metric and $\Omega$ the globally defined nondegenerate holomorphic $n-form$).\\
1) We say that a submanifold $L\subset X$ is {\em Special Lagrangian} if it is Lagrangian (of maximal dimension) with respect to $\omega$, and there exists a complex number of the form $e^{i\theta}$ 
such that $Im(e^{i\theta}\Omega)|_{L} = 0$. Such a $\theta$ is called the {\em phase} of the special lagrangian submanifold.\\
2) We say that a smooth map $f:X\to B$ to a smooth manifold $B$ of (real) dimension $n$ is a {\em Special Lagrangian Fibration} if $f$ is a submersion and for all $q\in B$ the submanifold $L_q = f^{-1}(q)\subset X$ is a special lagrangian submanifold of  
$(X,\omega,\mathbf{g},\Omega)$. We require also that the phase of the fibres is constant. 
\end{dfn}
\begin{lem}
Let $(X,\omega_X,\mathbf{g}_X,\Omega)_X$ be  a Calabi-Yau  manifold of  complex dimension $n$
Let $f:X\to B$ be a special lagrangian fibration with compact connected fibres, such that the metric of $X$ restricted to any fibre is flat. Then  $f$ is Riemannian.
\end{lem}
\textit{Proof}\\
In view of the description of deformations of special lagrangian manifolds of  ~\cite{ML}, it is enough to observe that harmonic forms on a flat manifold are covariant constant, and also their dual vector fields are covariant constant. As parallel transport is an isometry on any Riemannian manifold, and the complex involution is also an isometry, this implies that on each fibre you have an orthonormal frame of vector fields, whose transformations under the complex involution give a complete set of first order normal deformations of the fibre itself. This clearly implies that $f$ is Riemannian.
\ep

\section{$s$-\ka Manifolds}
\label{sec:ska}
In this section we introduce s-\ka manifolds. Contrary to almost s-\ka manifolds, they are extremely "rigid" objects, and it is more difficult to build examples. They enjoy however an extremely rich set of properties, especially when they are compact. They might be though of as "maximally symmetric" almost s-\ka manifolds.
\begin{dfn}
Let $M$ be a smooth manifold. We will indicate with 
$(~,~)_{\mathbf{g}}$ or with $(~,~)$ without further specification, 
the scalar product induced by a metric $\mathbf{g}$ on the tangent 
space of $M$, on its dual and on all their tensor powers (exterior, 
symmetric, etc.). We will also use the notation $|~|$ for the related 
(pointwise) norm, i.e. $|\alpha|~=~\sqrt{(\alpha,\alpha)}$.
\end{dfn}
\begin{dfn}[Definition of s-\ka manifold]
\label{dfn-ska}
A smooth manifold $M$ of dimension $n(s+1)$ together with a 
Riemannian metric $\mathbf{g}$ and $2$-forms 
$\omega_{1},...,\omega_{s}$ is {\em s-\ka} if the data satisfies the 
following property:
For each point of $M$ there exist an open neighborhood $\mathcal{U}$ 
of $p$ and a system of coordinates 
$x_{i},y^{j}_{i}$,$i~=~1,...,n$, $j~=~1,...,s$
on $\mathcal{U}$ such that:\\
1) $\forall j ~ ~\omega_{j}~=~\sum_{i}dx_{i}\wedge dy^{j}_{i}$,\\
2) $\mathbf{g}_{(\mathbf{x},\mathbf{y})}~=~\sum_{i}dx_{i}\otimes 
dx_{i}~+~\sum_{i,j}dy^{j}_{i}\otimes d 
y^{j}_{i}~+~\mathbf{O}(2)$.\\
Any such system of coordinates is called {\em standard}(s-\ka).
\end{dfn}
Note that the forms $\omega_i$ of an s-\ka manifold are closed, 
because they are constant in any standard coordinate system. 
Actually, more is true:
\begin{rmk}
If $(M,\omega_{1},...,\omega_{s},\mathbf{g})$ is an s-\ka manifold, 
$(M,\omega_{1},...,\omega_{s})$ is a (non degenerate) polysymplectic 
manifold 
\end{rmk}
The notion of almost $s$-\ka manifold, which we considered in the previous section, is a weakening of the notion of s-\ka manifold. To recover the full strength of the definition of $s$-\ka manifold, one 
then needs some integrability condition. 
\begin{rmk}
An s-\ka manifold is almost s-\ka
\end{rmk}
\begin{teo}
\label{almost-ka}
Let $M$ be a smooth manifold of dimension $n(s+1)$, let 
$\omega_{1},...,\omega_{s}$ be $2$-forms on $M$, and let $\mathbf{g}$ 
be a Riemannian metric on $M$. Then the following facts are 
equivalent:\\
$1)$~$(M,\omega_{1},...,\omega_{s},\mathbf{g})$ is an s-\ka 
manifold.\\
$2)$~$(M,\omega_{1},...,\omega_{s},\mathbf{g})$ is an almost s-\ka 
manifold, and
$\forall j~ ~\nabla\omega_{j}~=~0$, where $\nabla$ is the Levi-Civita 
(i.e. torsion free) connection associated to the metric $\mathbf{g}$.
\end{teo}
\textit{Proof} 
$1)\Rightarrow 2)$: 
This implication is an easy consequence of the fact that the 
Christoffel symbols (of the Levi-Civita connection) relative to any 
local coordinate system contain only the first derivatives of the 
metric.\\
$2)\Rightarrow 1)$: We will consider separately the cases $s = 1$ and  $s 
> 1$;\\
$s = 1$;\\ 
This case is essentially classical. We provide a proof for lack of a reference.
Let $(M,\omega,\mathbf{g})$ be a smooth manifold with a symplectic 
structure $\omega$ and a Riemannian metric $\mathbf{g}$, such that 
$\nabla\omega = 0$, where $\nabla$ is the Levi-Civita connection 
associated to $\mathbf{g}$, and such that the form $\omega$ and the 
metric $\mathbf{g}$ are compatible on $T_pM$ for all points $p\in M$. 
It is easy to see that the tensor $J$ defined as $\mathbf{g}(x,y) = 
\omega(x,Jy)$ is an almost complex structure on $M$, that 
$\mathbf{g}$ is compatible with $J$, and that $\omega$ is the 
fundamental form associated to the metric (see ~\cite{KN}, Volume II, Page 147). 
Therefore, $(M,J,\mathbf{g})$ is an almost Hermitian manifold. From 
the fact that $\nabla\omega = 0$ (and $\nabla\mathbf{g} = 0$, as is 
always the case) we deduce that $\nabla J = 0$, i.e. the Levi-Civita 
connection is almost complex. Using (~\cite{KN}, Volume II, Theorem 4.3 Page 148) 
we deduce that the almost complex structure is 
integrable, and moreover the metric is \ka with respect to it, as the 
fundamental 2-form $\omega$ is closed. Therefore, we 
can conclude that $(M,J,\mathbf{g})$ is a \ka manifold, with 
associated fundamental form $\omega$. Pick then complex coordinates 
$x_1+iy_1,...,x_n+iy_n$, such that in these coordinates 
\[\mathbf{g} = \sum_{ij}h_{ij}\left(dx_i\otimes dx_j + dy_i\otimes 
dy_j\right),~ ~\omega = \sum_{ij}h_{ij}dx_i\wedge dy_j\]
with~ $h_{ij} = \delta_{ij} + O(2)$. Their existence is guaranteed by 
the  \ka property. Moreover, 
\[0 = d\omega = \sum_{ijk}\frac{\partial h_{ij}}{\partial 
x_k}dx_k\wedge dx_i\wedge dy_j ~+~ \frac{\partial h_{ij}}{\partial 
y_k}dy_k\wedge dx_i\wedge dy_j\]
from which it follows that $\frac{\partial h_{ij}}{\partial x_k} = 
\frac{\partial h_{kj}}{\partial x_i}$ and $\frac{\partial h_{ij}}{\partial 
y_k} = 
\frac{\partial h_{kj}}{\partial y_i}$ for all $i,j,k \in\{1,...,n\}$.
Consider now the system of PDE's
\[\left\{\begin{array}{l}
\frac{\partial \tilde{x}_j}{\partial x_i} = h_{ij},~ ~i,j = 1,...,n
\\
~ ~\\
\frac{\partial \tilde{x}_j}{\partial y_i}(0) = \frac{\partial^2 
\tilde{x}_j}{\partial y_i y)k}(0) = 0,~ ~i,j,k = 1,...,n\\
~ ~\\
\tilde{x}_j(0) = 0,~ ~ j = 1,...,n
\end{array}\right.\]
The conditions $\frac{\partial h_{ij}}{\partial x_k} = \frac{\partial 
h_{kj}}{\partial x_i}$ and $\frac{\partial h_{ij}}{\partial 
y_k} = 
\frac{\partial h_{kj}}{\partial y_i}$ guarantee that the functions
\[\tilde{x}_j = \int_0^{x_1}h_{1j}(t,0,..,0)dt  + \cdots + 
\int_0^{x_n}h_{nj}(x_1,x_2,..,x_{n-1},t)dt,~ ~j=1,...,n\]
are solutions to the system above. Therefore, we can use 
$\tilde{x}_1,..,\tilde{x}_n,y_1,...,y_n$ as coordinates in a 
neighborhood of the point corresponding to $x_1=\cdots = y_n = 0$. 
Let $k_{ij}$ be the inverse matrix to $h_{ij}$. We then have that 
$d\tilde{x}_j = \sum_i h_{ij}dx_i + \sum_{i}\frac{\partial \tilde{x}_{j}}{\partial y_i}dy_{i}$. 
It follows that with these coordinates
\[\mathbf{g} = \sum_{ij}\left(k_{ij}d\tilde{x}_i\otimes d\tilde{x}_j - 
\sum_{k}k_{ij}\frac{\partial \tilde{x}_{j}}{\partial 
y_k}dy_{k}\otimes d\tilde{x}_{i} +
h_{ij}dy_i\otimes dy_j\right),~ ~ ~\omega = \sum_j d\tilde{x}_j\wedge 
dy_j\]
The expression for $\omega$ is obtained using the equations $\frac{\partial \tilde{x}_j}{\partial y_i} = \frac{\partial 
\tilde{x}_i}{\partial x_j}$. It remains to be checked that $\mathbf{g} = \Delta ~+~O(2)$ (where 
$\Delta$ is the identity matrix).  First, observe that $(Xk_{ij})_0 = 0$ for any vector field $X$, as we know 
that $h_{ij} = \delta_{ij} + O(2)$. Moreover, by construction all first order derivatives of $\frac{\partial \tilde{x}_{j}}{\partial 
y_k}$ vanish at the origin, and therefore the thesis follows.\\
$s > 2$;\\
Let now $M$ be a smooth manifold of dimension $n(s+1)$, with $s > 1$, 
and  let $\omega_{1},...,\omega_{s}$ and $\mathbf{g}$ be as defined in 
condition $2)$. Let $p$ be a point of $M$. Pick any standard 
polysymplectic coordinate system $x_{i},y^{j}_{i}$,$i~=~1,...,n$, 
$j~=~1,...,s$ centered at $p$, defined
on a neighborhood $\mathcal{U}$ of $p$ and such that:\\
1) $\forall j ~ ~\omega_{j}~=~\sum_{i}dx_{i}\wedge dy^{j}_{i}$,\\
2) $\mathbf{g}_{p} = \sum_i dx_idx_i~+~\sum_{ij}dy_i^jdy_i^j$\\
i.e. such that the induced coframe on $T_pM$ is orthonormal. Such a 
coordinate system exists from the definition of almost $s$-\ka 
manifold and from Theorem ~\ref{teo-multisympnorm}.
From the fact that $\nabla\omega_j = 0$ for all $j$, we deduce that 
parallel transport preserves the polysymplectic structure, and 
therefore it must preserve also the standard subspaces associated to 
it, among  which are the
\[<\frac{\partial}{\partial 
y^{1}_{1}},\ldots,\frac{\partial}{\partial y^{1}_{n}}>,\ldots,<\frac{\partial}{\partial 
y^{s}_{1}},\ldots,\frac{\partial}{\partial y^{s}_{1}}>\] 
From this we deduce that for any vector field $X$
\[\nabla_{X}\frac{\partial}{\partial 
y^{1}_{i}} = \sum_{l}dy^{1}_{l}\left(\nabla_{X}\frac{\partial}{\partial 
y^{1}_{i}}\right)\frac{\partial}{\partial 
y^{1}_{l}},\ldots,
\nabla_{X}\frac{\partial}{\partial 
y^{s}_{i}} = \sum_{l}dy^{s}_{l}\left(\nabla_{X}\frac{\partial}{\partial 
y^{s}_{i}}\right)\frac{\partial}{\partial 
y^{s}_{l}}\]
As a consequence, 
$\nabla_{\frac{\partial}{\partial x_i}}dx_l = -\sum_m 
\Gamma_{im}^ldx_m$, where $\Gamma_{im}^l = 
dx_l\left(\nabla_{\frac{\partial}{\partial 
x_i}}\frac{\partial}{\partial x_m}\right)$ are the usual Christoffel 
symbols. We will use the index notation $1,\ldots,n,$ $(11),\ldots,(ns)$ 
to indicate the $n(s+1)$ indices for the coordinates $x_{i},y^{j}_{i}$,$i~=~1,...,n$, 
$j~=~1,...,s$. The above considerations then amount to the fact that 
$\Gamma^{(ij)}_{\alpha m} = 0$ for any index $\alpha$, any numbers 
$i,m$ in the set $\{1,\ldots,n\}$ and any number $j$ in the set $\{1,\ldots,s\}$.
Consider now a coordinate change of the form
\[\tilde{x}_i = x_i + \sum_{mp}b^i_{mp}x_mx_p,~ ~ ~\tilde{y}^j_i = 
y^j_i(x_1,...,x_n,y^j_1,...,y^j_n)\]
where the functions $\tilde{y}^j_i$ are determined according to the 
Theorem of Caratheorody-Jacobi-Lie (~\cite{LM} Theorem 13.4 Page 136), so that $\omega_j = 
\sum_i d\tilde{x}_i\wedge d\tilde{y}^j_i$, and $\tilde{y}^j_i(0,...,0) 
= 0$. Note that it is crucial that the functions 
$\tilde{x}_{i}$ are in involution with respect to the Poisson 
structures associated (in the respective 
$x_{1},\ldots,x_{n},y^{j}_{1},\ldots,y^{j}_{n}$ spaces) to the various 
symplectic forms $\omega_{1},\ldots,\omega_{s}$. In view of the 
previous considerations, we see that also in the new coordinates we 
have $\nabla_{\frac{\partial}{\partial 
\tilde{x}_m}}\frac{\partial}{\partial 
\tilde{x}_p} = \sum_p\tilde{\Gamma}^l_{mp}\frac{\partial}{\partial 
\tilde{x}_l}$, if the $\tilde{\Gamma}$ are the Christoffel symbols 
in the new coordinates, and moreover\\
$\nabla_{\frac{\partial}{\partial 
\tilde{x}_m}}d\tilde{x}_l = 
\nabla_{\frac{\partial}{\partial 
x_m}}d\tilde{x}_l ~+~O(1) =$\\
 $\nabla_{\frac{\partial}{\partial 
x_m}}\left(dx_l + \sum_{ip}b^l_{ip}x_idx_p\right) ~+~O(1) =
 \sum_p\left(-\Gamma^l_{mp} + 
b^l_{mp}\right)d\tilde{x}_p ~+~O(1)$.\\ 
As it is also the case that $\nabla_{\frac{\partial}{\partial 
\tilde{x}_m}}d\tilde{x}_l = -\sum_p\tilde{\Gamma}^l_{mp}d\tilde{x}_p$,
if we choose $b^l_{mp}
= \Gamma^l_{mp}(0)$ (which we can do as the connection is torsion-free), 
we see that the symbols $\tilde{\Gamma}^l_{mp}$ in the new coordinate 
system vanish at the origin. For simplicity, we will indicate 
the new coordinates with $x_{i},y^{j}_{i}$, and the Christoffel 
symbols associated to them with $\Gamma$, dropping the tilde 
everywhere. 
We know also that for any index $\alpha$, and indicating with $(~ ~)_0$ the evaluation of a form at $0$,\\
$0 = \left(\nabla_{\alpha}\omega_{j}\right)_{0} = \left(\nabla_{\alpha}\sum_{i}dx_{i}\wedge 
dy^{j}_{i}\right)_{0} =  
\sum_{i}\left(dx_{i}\wedge (\nabla_{\alpha}dy^{j}_{i})\right)_{0} =$.\\
$-\sum_{i}\left(dx_{i}\wedge (\sum_{m k}\Gamma^{(ij)}_{\alpha (m 
k)}(0)dy^{k}_{m}~+~
\sum_{m}\Gamma^{(ij)}_{\alpha m}(0)dx_{m})\right)_{0}$\\
From this we deduce that $\Gamma^{(ij)}_{\alpha (m 
k)}(0) = 0$  and $\Gamma^{(ij)}_{\alpha m}(0) = \Gamma^{(mj)}_{\alpha 
i}(0)$ for all $i,j,k,m,\alpha$. We consider therefore the change of 
coordinates
\[\tilde{y}^{j}_i = y^{j}_i + \sum_{mp}\Gamma^{(ij)}_{mp}(0)x_mx_p,~ ~ 
~\tilde{x}_i = x_{i}\]
In the new coordinates we have
\[\sum_{i}dx_{i}\wedge d\tilde{y}^{j}_{i} = \sum_{i}dx_{i}
\wedge (dy^{j}_i + \sum_{mp}\Gamma^{(ij)}_{mp}(0)x_mdx_p) = 
\omega_{j},\]
as we showed before that $\Gamma^{(ij)}_{mp}(0) = \Gamma^{(pj)}_{mi}(0)$. All the equations for the 
Christoffel symbols that we have deduced so far still hold, because 
we did not make any assumption on the $y^{j}_{i}$ when we obtained 
them, apart from the fact that we were in polysymplectic coordinates.
Moreover, we have that 
\[\left(\nabla_{\frac{\partial}{\partial x_l}}d\tilde{y}^{j}_{i}\right)_{0} = 
\left(\nabla_{\frac{\partial}{\partial x_l}}dy^{j}_i + 
\sum_{mp}\Gamma^{(ij)}_{mp}(0)\nabla_{\frac{\partial}{\partial x_l}}(x_mdx_p)\right)_{0} = \]
\[\left(-\sum_{m}\Gamma^{(ij)}_{lm}(0)dx_{m} + 
\sum_{p}\Gamma^{(ij)}_{lp}(0)dx_{p}\right)_{0} = 0\]
From the previous equation, the symmetry of the Christoffel symbols 
coming from the fact that the connection is torsion-free, and the 
vanishing properties proved above, we see that all the 
Christoffel symbols 
vanish at $0$.\\
We know from the compatibility of the polysymplectic structure with 
the metric that there is a linear change of coordinates which sends 
the given coframe at $0$ to an orthonormal (but still polysymplectic) 
one. It follows that the same linear change, applied to the 
functions $x_{i},y^{j}_{i}$ will preserve the polysymplectic 
property, and will make the  coframe at $0$ orthonormal. Moreover, will not 
disrupt 
the vanishing property (at $0$) of the Christoffel symbols.\\
On the other hand, from the vanishing at the origin of all the Christoffel symbols (and 
the fact that the coordinate coframe at $0$ is orthonormal) it is 
straightforward to deduce that 
$\mathbf{g}~=~\sum_{i}dx_{i}\otimes dx_{i}~+~\sum_{i,j}dy^{j}_{i}
\otimes dy^{j}_{i}~+~\mathbf{O}(2)$.
\ep
\begin{cor}
Let $M$ be a smooth manifold, let $\omega_{1},...,\omega_{s}$ be 
$2$-forms and $\mathbf{g}$ be a Riemannian metric on it. Assume that 
$s\not= 2$. Then if for all $p\in M$ the forms induce a non 
degenerate polysymplectic structure on $T_pM$, the metric 
$\mathbf{g}_p$ is compatible with it and $\nabla\omega_j = 0$ for all 
$j = 1,..,s$ (where $\nabla$ is the Levi-Civita connection associated 
to $\mathbf{g}$), we have that 
$(M,\omega_{1},...,\omega_{s},\mathbf{g})$ is s-\ka
\end{cor}
\textit{Proof} From the fact that $\nabla\omega_j = 0$ for all $j = 
1,..,s$ we deduce that $d\omega_j = 0$ for all $j = 1,..,s$. Then, 
using Theorem ~\ref{simpsdiff2}, we know that 
$(M,\omega_{1},...,\omega_{s})$ is polysymplectic, and therefore by 
definition $(M,\omega_{1},...,\omega_{s},\mathbf{g})$ is almost 
s-\ka. Using the previous theorem, we conclude that 
$(M,\omega_{1},...,\omega_{s},\mathbf{g})$ is s-\ka.
\ep
\begin{rmk}
Let $(M,\me)$ be a Riemannian manifold. Then $\st^{*}(M)$, 
together with its canonical almost s-\ka structure (described in Example ~\ref{cotan-almostska}), is s-\ka if and only if $\me$ is everywhere flat
\end{rmk}
The proof of this remark is a long but straightforward computation in local coordinates.
\section{Examples of s-\ka manifolds}
\label{sec:exeska}
As we mentioned before, s-\ka manifolds are very rigid, and it is not easy to build examples when $s >1$. For $s=1$, we have the following very classical fact:
\begin{teo}
\label{1ka-ka}
Let $(M,\omega,\mathbf{g})$ be a 1-\ka manifold. Then  the equation
\[\forall x~y~\in~\mathbf{T}M~ ~\omega(x,Jy)~=~\mathbf{g}(x,y)\]
defines a complex structure on $M$, for which $\mathbf{g}$ is a \ka 
metric.\\
Conversely, if $(M,J,\mathbf{g})$ is a \ka manifold, the equation
\[\forall x~y~\in~\mathbf{T}M~ ~\omega(x,Jy)~=~\mathbf{g}(x,y)\]
defines a two-form on $M$, such that $(M,\omega,\mathbf{g})$ is a 
1-\ka manifold.
\end{teo}
\textit{Proof}
Let $(M,\omega,\mathbf{g})$ be a 1-\ka manifold. The equation
$\forall x~y~\in~\mathbf{T}M$ 
$\omega(x,Jy)~=~\mathbf{g}(x,y)$
defines a complex structure $J$ on $M$. Indeed,
it is clear from the non degeneracy and smoothness of both $\omega$ 
and $\mathbf{g}$ that the equation above defines a smooth section of 
$\mathbf{T}M^{*}\otimes\mathbf{T}M$, which induces an isomorphism on 
$TM_{p}$ for all $p~\in~M$. We first show that this $J$ is an almost 
complex structure on $M$. As this is a pointwise statement, we can 
choose a point $p$ on $M$ and standard coordinates centered at $p$ 
such that if the associated coframe is $dx_{i},dy^{j}_{i}$, on 
$\mathbf{T}M$ (over the coordinate neighborhood) we have
$\forall j ~\omega_{j}|_{p}~=~\sum_{i}dx_{i}\wedge dy^{j}_{i}$, and 
$\mathbf{g}_{p}~=~\sum_{i}dx_{i}\otimes 
dx_{i}~+~\sum_{i,j}dy^{j}_{i}\otimes dy^{j}_{i}~+~\mathbf{O}(2)$.
Therefore, using the notation
\[J^{0}~=~\sum_{i}\left(dx_{i}\otimes\frac{\partial}{\partial 
y^{1}_{i}}~-~
dy_{i}^{1}\otimes\frac{\partial}{\partial x_{i}}~+~
dy_{i}^{2}\otimes\frac{\partial}{\partial y^{3}_{i}}~-~
dy_{i}^{3}\otimes\frac{\partial}{\partial y^{2}_{i}}\right),\]
we must have that $J-J^{0}~=~\mathbf{O}(2)$.
This shows that $J^{2}~=~Id$, as this is clearly true over $p$, and 
the choice of $p$ was arbitrary. To show the integrability of $J$, it 
is enough, from the Theorem of Newlander and Niremberg, to show that 
the torsion vanishes. This is however immediate, as in the expression 
of the torsion $\mathbf{N}(X,Y)$ applied to two vector fields only 
the first derivatives of $J$ appear, and $J$ osculates to degree two 
to $J^{0}$, for which the torsion is clearly zero.
To show that the metric $\mathbf{g}$ is \ka with respect to $J$, it 
is enough to show that $\nabla J~=~0$, where $\nabla$ is the 
Levi-Civita connection relative to $\mathbf{g}$. If we covariantly 
differentiate the equation defining $J$, we obtain
$0~=~\nabla\mathbf{g}~=~\nabla\left(\omega\circ (Id\otimes J)\right)$.
From this, the fact that on any 1-\ka manifold $\nabla\omega~=~0$, 
and the non degeneracy of $\omega$, we conclude immediately (using 
"the Leibnitz rule") that
$\nabla J~=~0$, as desired.\\
In the other direction, assume that $(M,J,\mathbf{g})$ is a \ka 
manifold.
The equation
$\forall x~y~\in~\mathbf{T}M$ 
$\omega(x,Jy)~=~\mathbf{g}(x,y)$
defines a $2$-form $\omega$ on $M$, usually called the fundamental 
$2$-form (cfr. ~\cite[Volume II,Page 147]{KN}). The fact that 
$\omega$ is closed is an alternative definition of the \ka condition  
(cfr. ~\cite[Volume II, Page 149]{KN}) and from the fact that $\nabla 
J=0$ (cfr. ~\cite[Volume II, Corollary 4.4 Page 149]{KN}) and 
$\nabla\mathbf{g}=0$ it follows that $\nabla\omega=0$. Finally, the 
non degeneracy of $\mathbf{g}$ implies that $\omega$ is a symplectic 
form. We can now conclude using Theorem ~\ref{almost-ka}.
\ep

We now show that $3$-\ka~implies \hka. We take as definition of a \hka ~
manifold
the one given in ~\cite[Bryant, Definition 7 Page 156]{FU}.
\begin{teo}
\label{3ka-hka}
Let $(M,\omega_{1},\omega_{2},\omega_{3},\mathbf{g})$ be a 3-\ka 
manifold. Then the equations:\\
$1)~\forall x,y~\in~TM~\forall~i~\in\{1,2,3\}~ 
~\mathbf{g}(J_{i}x,J_{i}y)~=~\mathbf{g}(x,y),~ ~
\omega_{i}(J_{i}x,J_{i}y)~=~\omega_{i}(x,y)$,\\
$2)\forall x,y~\in~TM$
~$\omega_{1}(x,J_{1}y)~+~\omega_{1}(J_{2}x,J_{3}y)~=~\mathbf{g}(x,y)$,\\
~$\omega_{2}(x,J_{2}y)~+~\omega_{1}(J_{3}x,J_{1}y)~=~\mathbf{g}(x,y)$ 
and 
~$\omega_{3}(x,J_{3}y)~+~\omega_{3}(J_{1}x,J_{2}y)~=~\mathbf{g}(x,y)$\\
$3)~J_{1}^{2}~=~J_{1}^{2}~=~J_{1}^{2}~=~-Id,~ ~J_{1}J_{2}~=~J_{3}$\\
define uniquely three smooth sections $J_1,J_2,J_3$ of $TM^{*}\otimes 
TM$, which are complex structures, and which together with the metric 
$\mathbf{g}$ determine on $M$ a \hka~ structure via the definitions 
$\Omega_i(x,y)~=~-\mathbf{g}(x,J_i(y))$.
\end{teo}
{\it Proof} Let $(M,\omega_{1},\omega_{2},\omega_{3},\mathbf{g})$ be 
a 3-\ka manifold. To check that equations $1)-3)$ admit a (unique) 
solution, it is enough to check it at the center $p\in M$ of a system 
of standard coordinates. If $v_i,u^i_j$ is the orthonormal coframe 
over $p$ associated to the standard coordinates, we have 
$(\omega_j)_p=\sum_i v_i\wedge u^i_j$. One can verify directly that 
the equations $1)-2)$ (restricted over $p$) admit the unique 
solution:\\
$J_1(v_i)=u_1^i,~ ~J_1(u_1^i)=-v_i,~ ~J_1(u_2^i)=u_3^i,~ 
~J_1(u_3^i)=-u_2^i$\\
$J_2(v_i)=u_2^i,~ ~J_2(u_2^i)=-v_i,~ ~J_2(u_3^i)=u_1^i,~ 
~J_2(u_1^i)=-u_3^i$\\
$J_3(v_i)=u_3^i,~ ~J_3(u_3^i)=-v_i,~ ~J_3(u_1^i)=u_2^i,~ 
~J_3(u_2^i)=-u_1^i$,\\
and that the $J_i$ thus defined satisfy also $3)$, and are smooth. 
Because in standard coordinates the metric osculates to degree $2$ to 
a flat one, we see that also the $J_i$ osculate to order $2$ to their 
flat analogues. From this, it is immediate to show that their 
torsions vanish, because in the expression of the torsion 
$\mathbf{N}_{J_i}(X,Y)$ applied to two vector fields only the first 
derivatives of $J_i$ appear. This shows, using ~\cite[Bryant, Theorem 
4, Page 156]{FU}, that the three tensors $J_i$ together with 
$\mathbf{g}$ induce on $M$ a \hka~ structure, via the definitions 
$\Omega_i(x,y)~=~-\mathbf{g}(x,J_i(y))$. Note that one can see directly from the construction that the forms $\Omega_j$ are constant with respect to the metric.
\ep
Note that if we only assume that we start with an {\em almost} $3$-\ka manifold, we still get three almost complex structures, which induce a quaternionic structure on every tangent space. The integrability condition of $3$-\ka manifolds is used then to prove that these almost complex structures are integrable.
\begin{exe}[Compact s-ka manifolds with arbitrary s]~\\
Consider the manifold $M=\rnum^{s+1}$, with metric given in standard coordinates $(x_0,...,x_s)$ as $\me_{ij} = \delta_{ij}r_j$. Let $\omega_j = \sqrt{r_0r_j}dx_0\wedge dx_j$. Then $(M,\omega_1,...,\omega_s,\me)$ and $(M/\znum^{s+1},\omega_1,...,\omega_s,\me)$ are s-\ka manifolds
\end{exe}
\textit{Proof} The forms determine a polysymplectic structure, and the metric is clearly compatible with it. It is therefore enough to show that $\nabla^\me\omega_j = 0$ for all $j\in\{1,...,s\}$. This last fact is clear the forms are closed, and the Christhoffel symbols for $\me$ vanish everywhere in standard coordinates.
\ep
\begin{exe} The examples of almost $s$-\ka manifolds built starting from riemannian special lagrangian fibrations of Calabi-Yau manifolds are actually $s$-\ka  if the special lagrangian fibrations are {\em covariant constant}. This is the content of Theorem ~\ref{specialtoska}, part $4$.
\end{exe}

\section{Lefschetz operators}
\label{sec:lefop}
In this section we define  a family of operators 
(together with their adjoints and associated commutators) which 
generalize to $s\geq 1$ the standard Lefschetz operator of \ka 
manifolds. Throughout the section, we assume fixed an s-\ka manifold $(M,\omega_1,...,\omega_s,\mathbf{g})$.
\begin{dfn}
\label{llambdah}
The operators $L_{i}$ and $\Lambda_{i}$ on $\Omega^{*}(M)$, for an 
s-\ka manifold $(M,\omega_{1},...,\omega_{s},\mathbf{g})$ are defined 
as
$L_{i}(\alpha)~=~\omega_{i}\wedge\alpha,~ ~\Lambda_{i}~=~L_{i}^{*},~ 
~H_i~=~[L_i,\Lambda_i]$.
\end{dfn}
\begin{rmk}
\label{rmk-adjoint}
To define the adjoint to $L$, we simply used the {\em pointwise} 
definition
$\forall \alpha,\beta~\in~\Omega^{*}(M)$~ ~
$(L_{i}\alpha,\beta)~=~(\alpha,\Lambda_{i}\beta)$. We did not assume 
that $M$ is compact or oriented.
Note that here the scalar product of two forms is a function on $M$, 
not a number.
\end{rmk}
We now prove a first group of identities, which can be used to show 
that we have a representation of the Lie algebra 
$\mathbf{sl}(s+1,\rnum)$ on the space of forms of an s-\ka manifold. In 
the next section we will complete these identities, in order to 
show that we have an analogous representation on the de Rham 
cohomology of an s-\ka manifold, at least when it is compact 
orientable (and this last representation is induced by the one on the 
forms).
\begin{dfn}
The operators $E_{i},~E_{i}^{j},I_{i},I_{i}^{j}$ on 
$\Omega^{*}(\mathcal{U})$, for a standard coordinate neighborhood 
$\mathcal{U}$ with standard coordinates $x_{i},y^{j}_{i}$, are 
defined as\\
$E_{i}(\alpha)~=~dx_{i}\wedge\alpha$,~ 
~$E^{j}_{i}(\alpha)~=~dy^{j}_{i}\wedge\alpha$,~ 
~$I_{i}~=~\frac{\partial}{\partial x_i}\contr\alpha$,~ 
~$I^{j}_{i}(\alpha)~=~\frac{\partial}{\partial 
y^{j}_{i}}\contr\alpha$.
\end{dfn}
The reasoning in the proofs that follow in this section is very similar to the one 
that applies to \ka manifolds, used for example in 
~\cite[Pages 106-114]{GH}.
\begin{lem}
\label{rel-among-e-i}
The following relations hold among the operators 
$\ei,\eik,\ii,\iik$:\\
1)  $ \ei\eik~=~-\eik\ei ,~ ~ \ii\iik~=~-~\iik\ii$ \\
2)  $ \ei\iik~=~-\iik\ei ,~ ~  \ii\eik~=~-\eik\ii $\\
3)  $ \ei I^{j}~=~-I^{j}\ei ,~ ~ for~i~\not=~j $\\
4)  $ \eik I_{l}^{m}~=~-I_{m}^{l}\eik ,~ ~ for~(i,k)~\not=~(l,m)$ \\
5)  $ \ei\ii~+~\ii\ei~=~Id ,~ ~ \eik\iik~+~\iik\eik~=~Id $
\end{lem}
{\it Proof}~
All these identities are easily verified, using the anti commutativity 
property of the wedge product.
\ep
\begin{lem}
\label{descriptionlihi}
1) $L_{i}~=~\sum_{i}E_{i}E^{j}_{i},~ 
~\Lambda_{i}~=~\sum_{i}I^{j}_{i}I_{i}~+~\mathbf{O}(2)$.\\
2) 
$[L_{k},\Lambda_{k}]~=~\sum_{i}\left(\eik\iik\ei\ii-\iik\eik\ii\ei\right)~+~\mathbf{O}(2)$\\
3) For $h~\not=~k$ 
$[L_{k},\Lambda_{h}]~=~\sum_{i}\left(\eik 
I^{i}_{h}\right)~+~\mathbf{O}(2)$\\
\end{lem}
\textit{Proof}~\\
1) is immediate from the definitions, and the fact that in standard 
coordinates the metric osculates to degree two to a flat one.\\
2) 
~$[L_{k},\Lambda_{k}]~=~\left(\sum_{i}\eik\ei\right)\left(\sum_{j}I^{j}I^{j}_{k}\right)~-~
\left(\sum_{j}I^{j}I^{j}_{k}\right)\left(\sum_{i}\eik\ei\right) + 
\mathbf{O}(2)~=$\\
$~=~\sum_{i\not=j}\left(\eik\ei 
I^{j}I^{j}_{k}~-~I^{j}I^{j}_{k}\eik\ei\right)~+~
\sum_{i}\left(\eik\ei\ii\iik~-~\ii\iik\eik\ei\right)~+ 
\mathbf{O}(2)~=$\\
$=~\sum_{i}\left(\eik\ei\ii\iik~-~\ii\iik\eik\ei\right)~+ 
\mathbf{O}(2)~=~
\sum_{i}\left(\eik\iik\ei\ii-\iik\eik\ii\ei\right) + 
\mathbf{O}(2)~$\\
3)
$[L_{k},\Lambda_{h}]~=~
\left(\sum_{i}\eik\ei\right)\left(\sum_{j}I^{j}I^{j}_{h}\right)~-~
\left(\sum_{j}I^{j}I^{j}_{h}\right)\left(\sum_{i}\eik\ei\right)~+ 
\mathbf{O}(2)~=$\\
$=~\sum_{i\not=j}\left(\eik\ei 
I^{j}I^{j}_{h}~-~I^{j}I^{j}_{h}\eik\ei\right)~+~
\sum_{i}\left(\eik\ei 
I^{i}I^{i}_{h}~-~I^{i}I^{i}_{h}\eik\ei\right)~+ 
\mathbf{O}(2)~=$\\
$=~\sum_{i}\left(\eik\ei 
I^{i}I^{i}_{h}~-~I^{i}I^{i}_{h}\eik\ei\right)~+ 
\mathbf{O}(2)~=~
\sum_{i}\eik I^{i}_{h}\left(\ei\ii+\ii\ei\right)~+ 
\mathbf{O}(2)~=$\\
$= ~\sum_{i}\left(\eik I^{i}_{h}\right)+ 
\mathbf{O}(2)~$
\ep
\begin{teo}
\label{skaid-part1}
~ \\
1) For any $k,h~$
$[L_{k},L_{h}]~=~0,~ ~[\Lambda_{k},\Lambda_{h}]~=~0$.\\
2) $[H_{k},L_{k}]~=~2L_{k},~
[H_{k},\Lambda_{k}]~=~-2\Lambda_{k}$\\
3) When $k~\not=~h~$,
$[H_{k},L_{h}]~=~L_{h}~$,
$[H_{k},\Lambda_{h}]~=~-\Lambda_{h}$\\
4) When $k~\not=~h~$,
$[[L_{k},\Lambda_{h}],L_{h}]~=~L_{k}~$,
$[[L_{k},\Lambda_{h}],\Lambda_{k}]~=~-\Lambda_{h}$\\
5) When $h,k,l$ are all different,
$[[L_{k},\Lambda_{h}],L_{l}]~=~0~$,
$[[L_{h},\Lambda_{k}],\Lambda_{l}]~=~0$
\end{teo}
\textit{Proof}~
In all the proofs for this theorem we will ignore the terms in 
$\mathbf{O}(2)$. This is correct because the relations that we want 
to prove involve operators which act pointwise, and therefore it is 
enough to prove them at the center of a standard coordinate system. 
We will also use without further reference the relations of the 
previous two lemmas. Moreover, in all the points the second equation can be obtained from the first by taking adjoints on both sides. We will consequently prove only the first equations.\\
1) 
Note that $L_{k}$ is a multiplication operator in the exterior 
algebra, by an element of degree 2, which is therefore central. This 
implies immediately that two such operators must commute with each 
other.\\
2) 
$[[L_{k},\Lambda_{k}],L_{k}]~=~
[\sum_{i}\left(\eik\iik\ei\ii-\iik\eik\ii\ei\right),L_{k}]~=$\\
$\left(\sum_{i}\left(\eik\iik\ei\ii-\iik\eik\ii\ei\right)\right)
\left(\sum_{j}E^{j}_{k}E^{j}\right)~-$\\
$~\left(\sum_{j}E^{j}_{k}E^{j}\right)
\left(\sum_{i}\left(\eik\iik\ei\ii-\iik\eik\ii\ei\right)\right)~=$\\
$\sum_{i\not=j}\left(\eik\iik\ei\ii E^{j}_{k}E^{j}-
\iik\eik\ii\ei E^{j}_{k}E^{j}-E^{j}_{k}E^{j}\eik\iik\ei\ii
+E^{j}_{k}E^{j}\iik\eik\ii\ei\right)~+$\\
$\sum_{i}\left(\eik\iik\ei\ii\eik\ei-
\iik\eik\ii\ei \eik\ei-\eik\ei\eik\iik\ei\ii
+\eik\ei\iik\eik\ii\ei\right)~=$\\
$\sum_{i}\left(\eik\iik\ei\ii\eik\ei-
\iik\eik\ii\ei \eik\ei-\eik\ei\eik\iik\ei\ii
+\eik\ei\iik\eik\ii\ei\right)~=$\\
$\sum_{i}\eik\ei\left(\iik\eik\ii\ei-\eik\iik\ei\ii\right)~+
~\sum_{i}\left(\eik\iik\ei\ii-\iik\eik\ii\ei\right)\eik\ei~=$\\
$\sum_{i}\left(\eik\ei+\eik\ei\right)~=~2L_{k}$\\
3) 
$[[L_{k},\Lambda_{k}],L_{k}]~=~
\sum_{i}\left(\eik\iik\ei\ii-\iik\eik\ii\ei\right),L_{h}]~=$\\
$\left(\sum_{i}\left(\eik\iik\ei\ii-\iik\eik\ii\ei\right)\right)
\left(\sum_{j}E^{j}_{h}E^{j}\right)~+$\\
$-~\left(\sum_{j}E^{j}_{h}E^{j}\right)
\left(\sum_{i}\left(\eik\iik\ei\ii-\iik\eik\ii\ei\right)\right)~=$\\
$\sum_{i\not=j}\left(\eik\iik\ei\ii E^{j}_{h}E^{j}-
\iik\eik\ii\ei E^{j}_{h}E^{j}-E^{j}_{h}E^{j}\eik\iik\ei\ii
+E^{j}_{h}E^{j}\iik\eik\ii\ei\right)~+$\\
$\sum_{i}\left(\eik\iik\ei\ii E^{i}_{h}\ei-
\iik\eik\ii\ei E^{i}_{h}\ei- E^{i}_{h}\ei\eik\iik\ei\ii
+ E^{i}_{h}\ei\iik\eik\ii\ei\right)~=$\\
$\sum_{i}E^{i}_{h}\ei\left(\eik\iik+\iik\eik\right)~=~L_{h}$\\
4) 
$[[L_{k},\Lambda_{h}],L_{h}]~=~\left(\sum_{i}\eik 
I^{i}_{h}\right)\left(\sum_{j}E^{j}_{h}E^{j}\right)~-~
\left(\sum_{j}E^{j}_{h}E^{j}\right)\left(\sum_{i}\eik 
I^{i}_{h}\right)~=$\\
$\sum_{i\not=j}\left(\eik 
I^{i}_{h}E^{j}_{h}E^{j}~-~E^{j}_{h}E^{j}\eik 
I^{i}_{h}\right)~+~\sum_{i}\left(\eik 
I^{i}_{h}E^{i}_{h}\ei~-~E^{i}_{h}\ei\eik 
I^{i}_{h}\right)~=$\\
$\sum_{i}\eik\ei\left(I^{i}_{h}E^{i}_{h}~+~E^{i}_{h}I^{i}_{h}\right)~=~L_{k}$\\
5) 
$[[L_{k},\Lambda_{h}],L_{l}]~=~\left(\sum_{i}\eik 
I^{i}_{h}\right)\left(\sum_{j}E^{j}_{l}E^{j}\right)~-~
\left(\sum_{j}E^{j}_{l}E^{j}\right)\left(\sum_{i}\eik 
I^{i}_{h}\right)~=$\\
$\sum_{i,j}\left(\eik 
I^{i}_{h}E^{j}_{l}E^{j}~-~E^{j}_{l}E^{j}\eik 
I^{i}_{h}\right)~=~0$
\ep
\begin{rmk}
The previous theorem holds unchanged if instead of an $s$-\ka manifold we use an almost $s$-\ka manifold. All the proofs remain unchanged, as for the purposes of the theorem we are interested only on pointwise properties.
\end{rmk}
\begin{cor}
\label{slonforms}
Given an almost  s-\ka manifold $(M,\omega_{1},...,\omega_{s},\mathbf{g})$  we 
induce a representation of the Lie algebra $\mathbf{sl}(s+1,\rnum)$ 
on 
$\Omega^{*}(M)$, via a map defined on the standard Chevalley basis 
$\mathbf{e}^i,\mathbf{f}^i,\mathbf{h}^i$ ($i=1,...,s$) of 
$\mathbf{sl}(s+1,\rnum)$ as:\\
$\mathbf{h}^i\to -H_i~+~H_{i+1}~for~i<s,~ ~\mathbf{h}^s\to -H_s$\\
$\mathbf{e}^i\to [L_i,\Lambda_{i+1}]~for~i<s,~ ~\mathbf{e}^s\to L_s$\\
$\mathbf{f}^i\to [\Lambda_{i},L_{i+1}]~for~i<s,~ ~\mathbf{f}^s\to 
-\Lambda_s$
\end{cor}
\textit{Proof} The verification of the relations on 
$\mathbf{e}^i,\mathbf{f}^i,\mathbf{h}^i$ dictated by the Cartan 
matrix of type $A_{s}$ is done by a straightforward application of 
Theorem ~\ref{skaid-part1}.
As an example, we compute some of them (we will write 
$\mathbf{e}^i,\mathbf{f}^i,\mathbf{h}^i$ for the images of these 
elements of $\mathbf{sl}(s+1,\rnum)$ in $End_\rnum(\Omega^{*}(M))$):\\
1) ($i<s$): 
$[\mathbf{e}^i,\mathbf{f}^i]~=~[[L_i,\Lambda_{i+1}],[\Lambda_{i},L_{i+1}]]~=$\\
$[L_{i+1},[\Lambda_{i},[L_i,\Lambda_{i+1}]]]+
[\Lambda_{i},[[L_i,\Lambda_{i+1}],L_{i+1}]]~=~[L_{i+1},\Lambda_{i+1}]+[\Lambda_{i},L_i]~=~\mathbf{h}^i$\\
2) $[\mathbf{e}^s,\mathbf{f}^s]~=~-[L_s,\Lambda_{s}]~=~\mathbf{h}^s$\\
3)($i<s$): 
$[\mathbf{e}^i,\mathbf{h}^i]~=~[H_i-H_{i+1},[L_i,\Lambda_{i+1}]]~=$\\
$-[L_i,[\Lambda_{i+1},H_i]]~-~[\Lambda_{i+1},[H_{i},L_i]]~+~[L_i,[\Lambda_{i+1},H_{i+1}]]~+~[\Lambda_{i+1},[H_{i+1},L_i]]~=$\\
$-[L_i,\Lambda_{i+1}]~+~2[L_i,\Lambda_{i+1}]~+~2[L_i,\Lambda_{i+1}]~-~[L_i,\Lambda_{i+1}]~=2\mathbf{e}^i$\\
4) $[\mathbf{e}^s,\mathbf{h}^s]~=~-[H_s,L_s]~=~2\mathbf{e}^s$\\
5)($i<s-1$): 
$[\mathbf{e}^{i+1},\mathbf{h}^i]~=~[H_i-H_{i+1},[L_{i+1},\Lambda_{i+2}]]~=$\\
$-[L_{i+1},[\Lambda_{i+2},H_i]]-[\Lambda_{i+2},[H_i,L_{i+1}]]+[L_{i+1},[\Lambda_{i+2},H_{i+1}]]+[\Lambda_{i+2},[H_{i+1},L_{i+1}]]~=$\\
$[L_{i+1},\Lambda_{i+2}]~-~[\Lambda_{i+2},L_{i+1}]~-~[L_{i+1},\Lambda_{i+2}]~+~
2[\Lambda_{i+2},L_{i+1}]~=~-\mathbf{e}^{i+1}$\\
6) 
$[\mathbf{e}^{s},\mathbf{h}^{s-1}]~=~[H_{s-1}-H_{s},L_{s}]~=~-\mathbf{e}^{s}$\\
We leave the remaining verifications to the reader.

\ep

\section{The action of $\mathbf{sl}(s+1,\rnum)$ on cohomology. Primitive forms}
\label{sec:hodgeska}
In this section we complete the set of  identities which we begun to 
describe in Theorem ~\ref{skaid-part1}. These last identities will 
allow us to show that we have a representation of the Lie algebra 
$\mathbf{sl}(s+1,\rnum)$ on the cohomology of an s-\ka manifold, 
induced by the representation on the space of forms described in 
Corollary ~\ref{slonforms}. This will be done showing that the 
Laplacian $\Delta_{d}$ commutes with the action of 
$\mathbf{sl}(s+1,\rnum)$. 
\begin{teo}[s-\ka identities]
\label{skaid}
 Let $(M,\omega_{1},...,\omega_{s},\mathbf{g})$ be a compact oriented 
s-\ka manifold. Then we have that:\\
1) $\forall k~ ~[L_k,d]~=~0$;\\
2)If we define 
$d^{c}_{k}~:=~[L_{k},d^{*}]$,
we have that
$\forall k~ ~dd^{c}_{k}~+~d^{c}_{k}d~=~0$;\\
3) $\forall k~ ~[L_{k},\Delta_{d}]~=~[\Lambda_{k},\Delta_{d}]~=~0$,
where $\Delta_{d}$ is the $d$-Laplacian relative to the metric 
$\mathbf{g}$ and to the orientation.
\end{teo}
{\it Proof}\\
1) This equation follows immediately from the fact that $d\omega_k = 
0$.\\
2)
If we write down the expression for $d^{c}_{k}$ in standard s-\ka coordinates centered 
at a point $p\in M$, we see that no derivative of the metric appears. Therefore, when we write down the expression   for $dd^{c}_{k}~+~d^{c}_{k}d$, only the first derivatives of the 
metric are involved. We skip the details, as they are completely analogous to those of, for example, ~\cite[Pages 111-115]{GH}.\\
It follows, as in the classical case of \ka manifolds, that  to prove the equation it is enough to 
reduce to the case of a constant metric. 
When the metric is flat, however, the 
equation is easily seen to be equivalent (using $1)$) to 
$[L_k,\Delta_d]=0$, which with a flat metric follows immediately from 
the fact that $\omega_k$ is constant in flat (orthonormal) 
coordinates.\\
3) The second equation is the adjoint of the first. The first one, 
once written down explicitely in terms of $d$ and $d^*$, follows 
immediately from points $1)-2)$.

\ep
\begin{cor}
Let $(M,\omega_{1},...,\omega_{s},\mathbf{g})$ be a compact s-\ka 
manifold. Then there is a canonical representation of the simple Lie 
algebra $\mathbf{sl}(s+1,\mathbf{C})$ on $H^*(M,\mathbf{C})$. 
\end{cor}
\begin{dfn} For a compact oriented s-\ka manifold $M$, we indicate 
with $\mathcal{H}^q(M,\mathbf{C})$ the complexification of the space 
of $\Delta_d$-harmonic forms. We indicate with 
$L_k,\Lambda_k,H_k,\mathbf{e}^k,\mathbf{f}^k,\mathbf{h}^k$ the 
operators of $\mathbf{sl}(s+1,\rnum)$ defined in the last part of 
section \ref{sec:lefop}, also when they are acting on 
$\mathcal{H}^q(M,\mathbf{C})$. In this context they will be thought of 
as elements of $\mathbf{sl}(s+1,\mathbf{C})$.
\end{dfn}
We will now describe a canonical decomposition of the space of forms. 
The linear spans in the following definitions are always intended 
over the complex numbers.
\begin{dfn}~\\
$\gs^{+}~=~<~[L_i,\Lambda_{i+j}]~|~i,j>0~>~+~<~L_{i}~|~i>0~>,$\\
$\gs^{-}~=~<~[L_{i+j},\Lambda_{i}]~|~i,j>0~>~+~<~\Lambda_{j}~|~j>0~>,$\\
$\ks~=~subalgebra~of~\gs~generated~by~<~[L_i,\Lambda_{i+j}]~|~i>0,j\not=0~>,$\\
$\ks^{+}~=~<~[L_i,\Lambda_{i+j}]~|~i,j~>0~>,~ 
~\ks^{-}~=~<~[L_{i+j},\Lambda_{i}]~|~i,j~>0~>$\\
$LW^{q}~=~v~\in~\mathcal{H}^q(M,\mathbf{C})~s.th.~\gs^{-}v~=~(0)$\\
$\ls~=~<~L_i~|~i=1,...,s~>$\\
\end{dfn}
\begin{dfn}
\label{dfn-primq}
$\primq~=~\mathcal{U}_{\ks}LW^{q}$
\end{dfn}
\begin{rmk}
\label{primq-in-wedgeq}
From  the definition it is clear that $\primq$ is a sub-$\ks$-module 
of $\mathcal{H}^q(M,\cnum)$
\end{rmk}
Because $\gs$ is split and simple, we know that any finite dimensional complex representation 
will be completely reducible, and that any irreducible finite 
dimensional module must me a highest weight module. We will actually 
use lowest weight modules, and this is clearly correct in the case 
of $\gs$.
From the definition, it is clear that $LW^{q}$ is just the space of 
lowest weight vectors contained inside $\mathcal{H}^q(M,\mathbf{C})$.
\begin{teo}
\label{lefschetz-like}
For all $q$ there is a natural decomposition of $\ks$-modules
\[\mathcal{H}^q(M,\cnum)~=~\bigoplus_{i\geq 0}S^i(l_{s})\prim^{q-2i}\]
\end{teo}
\begin{lem}
\label{pq-ks-sumodule}
1) The algebra $\ks$ is semisimple.\\
2) $\primq~=~{\mathcal{U}_{\ks^{+}}LW^{q}}$.\\
3) $\Lambda_{j}{\mathcal{U}_{\ks^{+}}LW^{q}}~=~(0)~\forall i,q$
\end{lem}
{\it Proof}\\
1) From the definition, and proceeding as we did for $\gs$, we see 
that 
actually $\ks~\cong~\mathbf{sl}(s,\mathbf{C})$. Note that if $s=1$ we 
have 
$\ks~=~(0)$\\
2) From the Poincare'-Birkhoff-Witt theorem applied to the lie 
algebra 
$\ks$, ordering in a way to put $\ks^{-}$ before $\ks^{+}$, we get 
immediately what we want. \\
3) From the Poincar\'{e}-Birkhoff-Witt theorem applied to the lie 
algebra 
$\gs$, ordering in a way to put 
$<\Lambda_{i}|i=1,\ldots,s>~<~\ks^{-}~<~\ks^{+}~<~<L_{i}|i=1,\ldots,s>$
and observing that 
$\Lambda_{j}{\mathcal{U}_{\ks^{+}}LW^{q}}~\subset~\mathcal{H}^{q-2}(M,\cnum)$
we get immediately what we want.
\ep

\textit{Proof of the theorem} \\
1) 
${\mathcal{U}_{\ks^{+}}LW^{q}}~\cap~\ls 
\mathcal{H}^{q-2}(M,\cnum)~=~(0)$.\\ 
$LW^{q}~\cap~\ls \mathcal{H}^{q-2}(M,\cnum)~=~(0)$, because if 
$v~\in~LW^{q}~\cap~\ls \mathcal{H}^{q-2}(M,\cnum)$, by 
decomposing $\gs \mathcal{H}^{q-2}(M,\cnum)$ into irreducible modules 
we may 
suppose that $v~=~\sum_{i}v_{i}$, with 
$v_{i}~\in~LW^{q}~\cap~\gs^{+}u_{i}$
and $u_{i}~\in~LW^{q_{i}}$ with $q_{i}~\leq~q-2$. But then we must 
have $v_{i}~=~0$ for all $i$, because a lowest weight module can 
contain only one nonzero lowest weight vector (which is a generator).
Now, it is clear form the Poincare'-Birkhoff-Witt theorem that $\ls 
\mathcal{H}^{q-2}(M,\cnum)$ is a $\ks$-module (it is generated by 
lowest 
weight vectors of degree smaller than or equal to  $q-2$). Therefore,
${\mathcal{U}_{\ks^{+}}LW^{q}}~\cap~\ls \mathcal{H}^{q-2}(M,\cnum)$
is a finite dimensional $\ks$-module. As $\ks$ is semisimple, it 
follows that this module must be a direct sum of lowest weight 
modules. 
Lowest weight vectors $v$ for $\ks$ are defined by the property that
$\ks^{-}v~=~0$
However, we know that for all elements of 
${\mathcal{U}_{\ks^{+}}LW^{q}}$, and in particular 
for any such $v$, 
$\Lambda_{j}v~=~(0)~\forall i$
It follows that any lowest weight vector for $\ks$ in 
${\mathcal{U}_{\ks^{+}}LW^{q}}$ must be 
inside $LW^{q}$, and therefore if
${\mathcal{U}_{\ks^{+}}LW^{q}}~\cap~\ls 
\mathcal{H}^{q-2}(M,\cnum)~\not=~(0)$
then
$LW^{q}~\cap~\ls \mathcal{H}^{q-2}(M,\cnum)~\not=~(0)$.\\
2)
${\mathcal{U}_{\ks^{+}}LW^{q}}~+~\ls 
\mathcal{H}^{q-2}(M,\cnum)~=~\mathcal{H}^q(M,\cnum)$\\
Use the Poincare'-Birkhoff-Witt theorem applied to the lie algebra 
$\gs$, and putting $\ls$ after $\ks^{+}$, and observe that 
$\mathcal{H}^{*}(M,\cnum)$ is generated over $\gs^{+}$ by lowest 
weight vectors.\\
3) From points 1) and 2) we deduce that 
$\mathcal{H}^q(M,\cnum)~=~\primq~\oplus~\ls 
\mathcal{H}^{q-2}(M,\cnum)$. We now proceed inductively on $q$, the 
case $q\leq 2$ being taken care of by points 1) and 2). We  assume 
inductively that $\mathcal{H}^{q-2}(M,\cnum)~=~\bigoplus_{i\geq 
0}S^i(l_{s})\prim^{q-2i-2}$. To conclude, using points 1) and 2), it 
is enough to observe that, if $i\not= j$, 
\[S^i(l_{s})\mathbf{\prim}^{q-2i}\cap 
S^j(l_{s})\prim^{q-2j}~=~(0)\]
from standard Lie algebraic arguments, and that 
\[\ls 
S^i(l_{s})\prim^{q-2i-2}~=~S^{i+1}(l_{s})\prim^{q-2i-2}~=~S^{i+1}(l_{s})\prim^{q-2(i+1)}\]
\ep
\section{s-Lefschetz theorems}
From the previous section it follows that the canonical $\sl$ action 
on the forms of an s-\ka manifold induces an action on de Rham 
cohomology. Starting from this observation, and using the results of 
the previous sections, we will deduce a result similar to Lefschetz's 
theorem for \ka manifolds.
\begin{teo}
\label{teo-lie-independence}
The multiplication map
\[S^{r}\left(\ls\right)~\otimes~\primq(M)~\rightarrow~
\mathcal{H}^{q+2r}(M,\cnum)\]
is injective for $r~\leq~n~-~q$
\end{teo}
The proof will be given at the end of this section. The space of 
primitive forms was defined in the previous section. Unless otherwise 
stated, in what follows we will always assume that the base field is 
the complex numbers $\cnum$. We first need 
some lemmas.\\
In the following lemma we will use the multi index notation 
$dy^{\alpha}_{I_{\alpha}}~=~dy^{\alpha}_{i_{i}}\wedge\cdots\wedge
dy^{\alpha}_{i_{r}}$
if $I_{\alpha}~=~(i_{1},\ldots,i_{r})$, and $dx^J = 
dx_{j_1}\wedge\cdots\wedge dx_{j_t}$ if $J = (j_1,...,j_t)$. For an 
ordered multi index $I = (i_1,...,i_r)$, we also use the notation 
$|I|=r$. $I^o$ indicates the ordered complement to the multi index, 
and the intersection of two ordered multi indices is the ordered 
multi index having as entries those common to the two initial ones.
\begin{lem}
\label{hi-degree}
If $x_1,..,x_n,y^1_1,...,y^s_n$ are standard s-\ka coordinates around 
a point $p$,
\[H_{k}(dy^{1}_{I_{1}}\wedge\cdots\wedge dy^{s}_{I_{s}}\wedge 
dx^{J})~=~
(|I_{k}\cap J|-|I_{k}^{o}\cap J^{o}|)
dy^{1}_{I_{1}}\wedge\cdots\wedge dy^{s}_{I_{s}}\wedge dx^{J} + 
\mathbf{O}(2)\]
and therefore, if $|I_{k}|~=~d_{k}$, $|J|~=~n-t$
\[H_{k}(dy^{1}_{I_{1}}\wedge\cdots\wedge dy^{s}_{I_{s}}\wedge 
dx^{J})~=~
(d_{k}~-~t)
dy^{1}_{I_{1}}\wedge\cdots\wedge dy^{s}_{I_{s}}\wedge dx^{J} + 
\mathbf{O}(2)\]
\end{lem}
{\it Proof}
For the first statement, recall the explicit description of the 
action of $H_{k}$ given in Lemma ~\ref{descriptionlihi},
$H_{k}~=~[L_{k},\Lambda_{k}]~=~\sum_{i}\left(\eik\iik\ei\ii-\iik\eik\ii\ei\right)  
+ \mathbf{O}(2)$.
With this description the statement is clear.
For the second one, suppose we have the following decomposition into 
disjoint subsets $A,B,C,D$ of $\{1,\ldots,n\}$:
\[I_{k}~=~A\cup B,~J~=~A\cup C,~A\cup B\cup C\cup D~=~ 
\{1,\ldots,n\}\]
Then
$I_{k}^{o}~=~C\cup D,~J^{o}~=~B\cup D,~I_{k}\cap J~=~A,~I_{k}^{o}\cap 
J^{p}~=~D$.
We have therefore that
$|I_{k}\cap J|~=~|A|,~|I_{k}^{o}\cap 
J^{p}|~=~|D|$.
We also know that
$(|A|+|B|)~~+~(|A|~+~|C|)~=~d_{k}+(n-t)$,
$|A|+|B|+|C|+|D|~=~n$.
Therefore, by subtracting the second equation from the first one,
$|A|~-~|D|~=~d_{k}-t$.
We conclude by observing that also
$|I_{k}\cap J|-|I_{k}^{o}\cap J^{o}|~=~|A|-|D|$
as desired.
\ep
\begin{lem}
\label{xj-hi-degree}
If $v$ is $H_{i}$-homogeneous, say $H_{i}v~=~\lambda_{i}v$, for some 
fixed $i$, then
$H_{i}L_{j}v~=~\left(\lambda_{i}~+~1~+~\delta_{ij}\right)L_{j}v$
\end{lem}
{\it Proof}
It is enough to observe that
$H_{i}L_{j}v~=~L_{j}H_{i}v~+~[H_{i},L_{j}]v,$
and that from the structure equations for $\gs$
$[H_{i},L_{j}]~=~(1+\delta_{ij})L_{j}$.
\ep
\begin{lem}
\label{sl2-injectivity}
Let $X,Y,H$ be nonzero linear operators on a finite dimensional 
vector space $M$ such that 
$[X,Y]~=~H,~[H,X]~=~2X,~[H,Y]~=~-2Y$
(i.e. such that $<X,Y,H>$ is an $\mathbf{sl}(2)$). Then if 
$v~\in~M$ is a vector such that $H(v)~=~-\lambda v$, $\lambda~>~0$, 
we 
can conclude that 
$X^{\lambda}v~\not=~0$
\end{lem}
{\it Proof}\\
We decompose $M$ into irreducible 
$\mathbf{sl}(2)$-modules, which we may assume to be lowest weight 
modules,
$M~=~\bigoplus_{i}\left(<~X^{k}~|~k\geq 0~>u_{i}\right),~Yu_{i}~=~0$
Then if
$v~=~\sum_{i}v_{i},~v_{i}~\in~\left(<~X^{k}~|~k\geq 
0~>u_{i}~\setminus (0)\right),$
we must have that $H(v_{i})~=~-\lambda v_{i}$ for all $i$. It is 
therefore enough to prove the statement for $M$ a lowest weight 
module of the form
$M~=~<~X^{k}~|~k\geq 0~>u,~Yu~=~0$.
At this point the proof is straightforward.
\ep
\begin{lem}
\label{inj-monomials}
Let $\mathbf{m}$ be a nonzero monomial of degree $n-q$ in the $L$'s, 
and let 
$x~\in~\mathcal{H}^q(M,\cnum)$, with $q\leq n$. Then
$\mathbf{m}x~=~0~\Rightarrow~x~=~0$
\end{lem}
{\it Proof}
For $q=0$ this is trivial, so assume $Q>0$. Assume also $x~\not= 0$.\\
We may as well assume that $x$ is homogeneous with respect to the 
$H_{i}$, because the equation $\mathbf{m}x~=~0$ can be decomposed in 
its homogeneous 
parts.  Let $p\in M$ be a point where $x_p\not= 0$. It is clearly 
enough to prove that $(\mathbf{m}x)_p\not= 0$. For this, take a set 
of standard coordinates $x_1,...,x_n,y^1_1,...,y^s_n$ around $p$, and 
let $\Omega = <(dx_1)_p,...,(dx_n)_p >$. We can from now on assume 
without loss of generality that we are working in an s-\ka vector 
space, with a fixed standard basis. We will continue to use the 
notation $H_j,L_j,\Lambda_j$ for the operators acting on this vector 
space, which are obtained by restriction from the operators 
$H_j,L_j,\Lambda_j$ acting on $\mathcal{H}^q(M,\cnum)$.\\
We can assume that the number of elements from $\Omega$ in each 
"monomial" of $x_p$ is fixed, say $n-t$, because $m$ will increase 
that 
number in each monomial by $n-q$. Therefore from Lemma 
~\ref{hi-degree} we 
know that for each monomial $y$ of $x$, and hence for $x$ itself, we 
have
$H_{i}x~=~\left(d_{i}-t\right)x,~\sum_{j}d_{j}~<~t$.
Suppose now that 
$\mathbf{m}~=~L_{s}^{b_{s}}\cdots 
L_{1}^{b_{1}}$,$\sum_{i}b_{i}~=~n-q$.
Applying Lemma ~\ref{sl2-injectivity} once to the triple 
$L_{1},\Lambda_{1},H_{1}$ we see that $L_{1}^{b_{1}}x~\not= 0$, 
because
$H_{1}x~=~\left(d_{1}-t\right)x,~t-d_{1}~\geq~t-\sum_{i}d_{i}~~=~n-q~\geq~b_{1}$.
Moreover, for $i~\not=~1$, we have from Lemma ~\ref{xj-hi-degree}
\[H_{i}L_{1}^{b_{1}}x~=~\left(d_{i}-t+b_{1}\right)x\]
Proceeding inductively, we can therefore assume that 
$L_{k}^{b_{k}}\cdots L_{1}^{b_{1}}x~\not=0$,
and for $i~>~k$
$H_{i}L_{k}^{b_{k}}\cdots L_{1}^{b_{1}}x~=~\left(d_{i}-t+\sum_{j\leq 
k}b_{j}\right)x$.
If one observes that\\
$-\left(d_{i}-t+\sum_{j\leq 
k}b_{j}\right)~=~t-d_{i}-\sum_{j\leq 
k}b_{j}\geq~t-\sum_{l}d_{l}-\sum_{j\leq 
k}b_{j}~=$\\
$(n-q)-\sum_{j\leq 
k}b_{j}~\geq~b_{k+1}$,\\
then applying Lemma ~\ref{sl2-injectivity} to the triple 
$L_{k+1},\Lambda_{k+1},H_{k+1}$ we see that\\
$L_{k+1}^{b_{k+1}}\cdots L_{1}^{b_{1}}x~\not=0,$
and applying Lemma ~\ref{xj-hi-degree} we see also that 
for $i~>~k+1$
$H_{i}L_{k+1}^{b_{k+1}}\cdots 
L_{1}^{b_{1}}x~=~\left((d_{i}-t+\sum_{j\leq 
k}b_{j})~+~b_{k+1}\right)x~=$
$\left(d_{i}-t+\sum_{j\leq 
k+1}b_{j}\right)x$,
concluding the inductive step.
\ep
\begin{dfn}
We will indicate with $\llex$ ( resp. $\sllex$) the  degree 
lexicographic  (respect. strict degree lexicographic) ordering on 
the multi degrees of monomials in the $L_{i}$'s. We will sometimes 
extend this to an ordering on the monomials themselves.
\end{dfn}
\begin{lem}
\label{shifting-ms}
Let $\mathbf{m}$ be a monomial in $\ls$ with respect to the standard 
basis. Then inside $\mathcal{U}_{\gs}$ we have:\\
1)~ For any pair of different indices $l~>~m$ 
we have that if $\mathbf{m}~\not=~L_{s}^{deg(\mathbf{m})}$ (i.e. 
$\mathbf{m}$ is not $\llex$-maximal)
$[H_{l,m},\mathbf{m}]~\in~span_{\ring}\left\{\mathbf{n}~|~\mathbf{m}~\sllex~\mathbf{n}
~\right\}$
and if $\psi~\in~\ks^{-}$ is generic, $[\psi,\mathbf{m}]~\not=~0$.
We have also that for any positive $r$
$[\ks^{-},L_{s}^{r}]~=~(0)$.
If we started with $l~<~m$, we would have that if
$\mathbf{m}$ is not $\llex$-maximal,
$[H_{l,m},\mathbf{m}]~\in~span_{\ring}\left\{\mathbf{n}~|~\mathbf{n}~\sllex~\mathbf{m}
~\right\}$.\\
2)~ For any $\phi~\in~\mathcal{U}_{\ks^{-}}$ we have
$\phi\mathbf{m}~\in~span_{\ring}\left\{\mathbf{n}~|~\mathbf{m}~
\llex~\mathbf{n}~\right\}\mathcal{U}_{\ks^{-}}$
\end{lem}
{\it Proof}\\
We will use the notation $H_{l,m} = [L_l,\Lambda_m]$, so that 
$H_{l,l} = H_l$.
We know from  Theorem
 ~\ref{skaid-part1} that 
\[[H_{l,m},L_{i}]~=~\left\{
\begin{array}{ll}
0 & for~i~\not=~m\\
L_{l} & for~i~=~m
\end{array}\right.
\]
If we want to evaluate $[H_{l,m},L_{i}\mathbf{n}]$, we can use the 
standard relation
\[[H_{l,m},L_{i}\mathbf{n}]~=~H_{l,m}[L_{i},\mathbf{n}]~+~[H_{l,m},L_{i}]\mathbf{n}\]
and proceed inductively on $deg(\mathbf{m})$. It is clear that if 
$\psi~\in~\ks^{-}$ is not in the span of
$\left\{~H_{i,j}~|~L_{j}\not|~\mathbf{m}~\right\}$
we have $[\psi,\mathbf{m}]~\not=~0$, and the above space is not all 
of $\ks^{-}$ precisely when 
$\mathbf{m}~\not=~L_{s}^{deg(\mathbf{m})}$.
For the proof of the second part, we have only to use the first part 
and the associativity property of multiplication.
\ep
\begin{lem}
\label{lie-independence-mono}
Let $v~\in~LW^{q}$
be a homogeneous element (also with respect to the action), let 
$\left\{\eta_{1},\ldots,\eta_{t}\right\}~\subset~\mathcal{U}_{\ks^{+}},~
\left\{\mathbf{m}_{1},\ldots\mathbf{m}_{t}\right\}~\subset~S^{r}\left(\ls\right)$
(where the $\mathbf{m}_{i}$ are assumed to be different monomials 
with respect 
to the given basis) and assume that
$\mathbf{m}_{1}\eta_{1}v~+~\cdots~+~\mathbf{m}_{t}\eta_{t}v~=~0$.
Then if $r~\leq~n~-~q$ we have that
$\eta_{1}v~=~\cdots~=~\eta_{1}v~=~0$.
\end{lem}
{\it Proof}
Suppose by contradiction that
$\mathbf{m}_{1}\eta_{1}v~+~\cdots~+~\mathbf{m}_{t}\eta_{t}v~=~0$
with all the $\mathbf{m}_{i}$ different monomials with respect 
to the given basis, and all the $\eta_{1}v~\not=~0$. We want to 
proceed by (inverted) induction on the lexicographic ordering of the 
multi degree of the 
monomials appearing, using the fact that the case $t~=~1$ is known 
(it is the content of Lemma ~\ref{inj-monomials}).
Suppose therefore that we have $t~>~1$, and let $|\mathbf{m}|$ 
indicate the multi degree of the monomial $\mathbf{m}$. We may assume 
without loss of generality that
$\forall~i>1~|\mathbf{m}_{1}|~\sllex~|\mathbf{m}_{i}|$.
Now, it is a standard fact of the theory of lowest weight modules 
for $\ks~\cong~\mathbf{sl}(s)$ (or for $\mathbf{sl}(s+1,\cnum)$) 
that there must be 
$\phi~\in~\mathcal{U}_{\ks^{-}}$ with
$\phi\eta_{1}v~=~v$.
Then, multiplying the relation assumed to exist by $\phi$, and using 
the second part of
Lemma ~\ref{shifting-ms}, we get a new relation
$\mathbf{m}_{1}v~+~\mathbf{m}^{\prime}_{2}\eta_{2}^{\prime}v~+
\cdots~+~\mathbf{m}^{\prime}_{t^{\prime}}\eta_{t}^{\prime}v~=~0$ with 
$|\mathbf{m}_{1}|~\sllex~|\mathbf{m}^{\prime}_{i}|$ for all $i$.
From the first part of the same Lemma, we obtain that for generic 
$\psi~\in~\ks^{-}$ the relation
$\psi\left(\mathbf{m}_{1}v~+~\mathbf{m}^{\prime}_{2}\eta_{2}^{\prime}v~+
\cdots~+~\mathbf{m}^{\prime}_{t^{\prime}}\eta_{t}v\right)~=~0$
is of the form
\[\mathbf{m}^{\prime}_{2}\eta_{t}^{\prime\prime}v~+
\cdots~+~\mathbf{m}^{\prime}_{t^{\prime}}\eta_{t}^{\prime\prime}v~+
\sum_{\mathbf{m}\not=~\mathbf{m}^{\prime}_{2},\ldots\mathbf{m}
^{\prime}_{t^{\prime}},|\mathbf{m}_{1}|\sllex|\mathbf{m}|}\eta_{\mathbf{m}}v~=~0\]
with all the $\eta_{t}^{\prime\prime}v$ different from zero.\\
By inverse $\llex$-induction, we conclude that we can reduce to the 
case of $t~=~1$, and therefore we obtain a contradiction.
\ep

{\it Proof of theorem ~\ref{teo-lie-independence}}\\
Let $w_1,...,w_t$ be elements of $\primq(M)$, with $q\leq n$, let 
$\left\{\mathbf{m}_{1},\ldots\mathbf{m}_{t}\right\}~\subset~S^{r}\left(\ls\right)$
be a set of distinct monomials (with respect 
to the given basis $L_1,...,L_s$), with $r = n-q$, and assume that
\[\mathbf{m}_{1}w_{1}~+~\cdots~+~\mathbf{m}_{t}w_{t}~=~0\]
Let  $\mathcal{B}$ be a fixed basis of $LW^q$ formed by homogeneous  
vectors (with respect to the $H_j$). From Definition  
~\ref{dfn-primq} we know that we can write any $w_i$ as a sum of the 
form $\sum_j\eta^{i,j}v_{k_{i,j}}$, where the $\eta^{i,j}$ are 
different elements of $\mathcal{U}_{\mathbf{k}_s^+}$  and the 
$v_{k_{i,j}}$ are elements of $\mathcal{B}$, all different for fixed 
$i$. The relation becomes
\[\sum_{i,j}\mathbf{m}_{i}\eta^{i,j}v_{k_{i,j}} = 0\]
By putting if necessary some of the $\eta^{i,j}$ equal to zero, and 
reordering their second index, we may assume that all the expressions 
for the $w_i$ involve all the elements from $\mathcal{B}$, and that 
$\eta^{i,j}$ is the coefficient relative to $v_j$. The equation above 
then becomes
$\sum_j\left(\sum_i\mathbf{m}_{i}\eta^{i,j}\right)v_j = 0$.
Because different $v_k$ generate distinct (and disjoint) 
$\slc$-modules, we can decompose the dependence relation above 
according to the elements of $\mathcal{B}$ involved, and we obtain 
therefore that
\[\forall ~j~ ~\sum_i\mathbf{m}_{i}\eta^{i,j}v_j = 0\]
From Lemma ~\ref{lie-independence-mono} we conclude that all the 
$\eta^{i,j}$ must be zero, and therefore all the $w_i$ are zero, as 
desired.
\ep
\begin{rmk} For the canonical $2$-\ka structure on the tree dimensional Torus $\textbf{T}^3$, we have that $dim_\rnum(\prim^2(\textbf{T}^3)) = 1$ and $n = 1$, and therefore it is not true in general that for a compact $2$-\ka manifold $M$ one has $\primq(M) = 0$ for $q > n$. Another counterexample to the vanishing is given by the orientation class on any compact oriented $s$-\ka manifold $M$ of dimension greater than $4n$, e.g. any torus $\textbf{T}^{s+1}$ with its canonical $s$-\ka structure and with $s \geq 4$. Indeed, such a class cannot be in the image under application of the $L_j$ operators of any class of degree less than $dim(M) - 2n > 2n$.
\end{rmk}
\section{Conclusions}
\label{sec:conclusions}
In this paper we tried to introduce the notions of polysymplectic manifold and of (almost) $s$-\ka manifold in the most straightforward way possible. However, it should be noted that there is a very natural way to generalize further the above definitions, and that this generalization may prove useful. For the sake of completeness, we give it here, so that it doesn't clutter the rest of the paper with unnecessary generality:
\begin{dfn}
Let $\mathbf{G}$ be a hypergraph on the vertex set $\mathcal{V}$, with edge set $\mathcal{E}$. Assume that for all vertices $v\in\mathcal{V}$ it is given an ordering $<_v$ on $\mathcal{E}(v)$, the set of edges containing the vertex $v$. Let $n$ be a natural number. A smooth manifold $M$ together with smooth differential forms $\omega_v$ (one for each element of $\mathcal{V}$) is said to be $\mathbf{G}$-symplectic of rank $n$ if for any $p\in M$ there is an open neighborhood of $p$ in $M$, and coordinates $x_i^e, i = 1,...,n, e\in \mathcal{E}$ such that 
\[\omega_v = \sum_{i=1}^n dx_i^{e_1(v)}\wedge\cdots\wedge dx_i^{e_{k(v)}(v)}\]
where $e_1(v) <_v \cdots <_v e_{k(v)}(v)$ are all the edges containing the vertex $v$ in increasing order.
\end{dfn}
\begin{rmk}
Polysymplectic manifolds are obtained from the hypergraph (with orderings)
\[\mathbf{P}_s = \left\{\{v_1\}_1,...,\{v_s\}_s,\{v_1,...,v_s\}_{s+1}\right\},~ ~\forall j~\{v_1,...,v_s\}_{s+1} <_{v_j} \{v_j\}_j\]
Note that the markings on the edges are necessary to distinguish edges containing the same vertices. This happens for $s = 1$ (symplectic manifolds): $\mathbf{P}_1 = \left\{\{v_1\}_1,\{v_1\}_2\right\},~ ~\{v_1\}_2 <_{v_1} \{v_1\}_1$. The hypergraph 
\[\mathbf{M}_s = \left\{\{v_1\}_1,...,\{v_1\}_s\right\},~ \{v_1\}_1 <_{v_1}\cdots <_{v_1} \{v_1\}_s\] 
describes a space with a single closed $s$-form. We believe that this object has been studied in the literature, under the name {\em  $s$-symplectic manifold}, or possibly {\em multisymplectic manifold}. This generalization of symplectic manifolds is in a sense ``transversal'' to the one that we studied in the present paper. Note that $\mathbf{P}_1 = \mathbf{M}_2$.
\end{rmk}
We do not write down  the easy generalization of the notions of almost $s$-\ka and $s$-\ka manifold, to that of {\em (almost) $\mathbf{G}$-\ka  manifold}. The interested reader can do that by himself, mimiking the definition of (almost) $s$-\ka manifold. We did not check how much of the theory generalizes to this setting. This might be a good exercise!\\
Concerning the study of special lagrangian fibrations, in a forthcoming paper we will show that the language of (almost) $2$-\ka manifolds can be very effective in studying their role in the approach of Strominger, Yau and Zaslow to mirror symmetry.

~ \\
Michele Grassi\\
Dipartimento di Matematica\\
Universita' di Pisa\\
via Buonarroti, 2\\
56127 Pisa, Italy\\
e-mail:~\textbf{grassi@dm.unipi.it}

\end{document}